\documentclass[11pt, a4paper]{article}
\usepackage{amsmath,amssymb,latexsym,color,enumerate}
\usepackage{authblk}
\usepackage[mathscr]{eucal}
\usepackage[colorlinks,
            linkcolor=blue,
            anchorcolor=green,
            citecolor=magenta
           ]{hyperref}
\newtheorem{thm}{Theorem}[section]

\newtheorem{lemma}[thm]{Lemma}
\newtheorem{cor}[thm]{Corollary}

\newtheorem{obs}{Observation}

\newtheorem{example}{Example}[section]
\newtheorem{defin}{Definition}[section]

\newcommand{\proof}{{\it Proof.\quad}}
\newcommand{\qed}{\hfill\Box\medskip}

\allowdisplaybreaks
\usepackage{CJK}

\begin{document}
%\begin{CJK*}{GBK}{song}

\renewcommand{\baselinestretch}{1.3}
%%%%%%%%%%%%%%%%%%%%%%%%%%%%%%%%%%%%%%%%%%%%%%%%%%%%%%%%%%%%%%%%%%%%%%%%%%%%%%%%%%%%%%%%
%%%%%%%%%%%%%%%%%%%%%%%%%%%%%%%%%%%%%%%%%%%%%%%%%%%%%%%%%%%%%%%%%%%%%%%%%%%%%%%%%%%%%%%%

\title{\bf  Some intersection theorems for finite sets}

\author[1]{Mengyu Cao\thanks{E-mail: \texttt{caomengyu@mail.bnu.edu.cn}}}
\author[1]{Mei Lu\thanks{E-mail: \texttt{lumei@tsinghua.edu.cn}}}
\author[2]{Benjian Lv\thanks{Corresponding author. E-mail: \texttt{bjlv@bnu.edu.cn}}}
\author[2]{Kaishun Wang\thanks{E-mail: \texttt{wangks@bnu.edu.cn}}}
\affil[1]{\small Department of Mathematical Sciences, Tsinghua University, Beijing 100084, China}
\affil[2]{\small Laboratory of Mathematics and Complex Systems (Ministry of Education), School of Mathematical Sciences, Beijing Normal University, Beijing 100875, China}
\date{}
\maketitle
\begin{abstract}
Let $n$, $r$, $k_1,\ldots,k_r$ and $t$ be positive integers with $r\geq 2$, and $\mathcal{F}_i\ (1\leq i\leq r)$ a family of $k_i$-subsets of an $n$-set $V$.
The families $\mathcal{F}_1,\ \mathcal{F}_2,\ldots,\mathcal{F}_r$ are said to be $r$-cross $t$-intersecting if $|F_1\cap F_2\cap\cdots\cap F_r|\geq t$ for all $F_i\in\mathcal{F}_i\ (1\leq i\leq r),$ and said to be non-trivial if $|\cap_{1\leq i\leq r}\cap_{F\in\mathcal{F}_i}F|<t$. If the $r$-cross $t$-intersecting families $\mathcal{F}_1,\ldots,\mathcal{F}_r$ satisfy $\mathcal{F}_1=\cdots=\mathcal{F}_r=\mathcal{F}$, then $\mathcal{F}$ is well known as $r$-wise $t$-intersecting family. In this paper, we describe the structure of non-trivial $r$-wise $t$-intersecting families with maximum size, and give a stability result for these families. We also determine the structure of non-trivial $2$-cross $t$-intersecting families with maximum product of their sizes.

\medskip

\noindent {\em AMS classification:}  05D05, 05A10

\noindent {\em Key words:} non-trivial $r$-wise $t$-intersecting family; non-trivial cross $t$-intersecting family

\end{abstract}

\section{Introduction}
Let $n$ and $k$ be integers with $1\leq k\leq n.$ Write $[n]=\{1,2,\ldots,n\}$ and denote by ${[n]\choose k}$ the family of all $k$-subsets of $[n].$ For a positive integer $t$, a family $\mathcal{F}\subseteq {[n]\choose k}$ is said to be \emph{$t$-intersecting} if $|A \cap B|\geq t$ for all $A, B\in\mathcal{F}.$ A $t$-intersecting family is called \emph{trivial} if all its members contain a common specified $t$-subset of $[n]$, and \emph{non-trivial} otherwise.

The famous Erd\H{o}s-Ko-Rado Theorem \cite{Erdos-Ko-Rado-1961-313} showed that each extremal $t$-intersecting family of ${[n]\choose k}$ consists of all $k$-subsets that contain a fixed $t$-subset of $[n]$ for $n>n_0(k,t)$. It is known that the smallest possible such function $n_0(k, t)$ is $(t+1)(k-t+1).$ This was proved by Frankl \cite{Frankl-1978} for $t\geq 15$ and subsequently determined by Wilson \cite{Wilson-1984} for all $t$. In \cite{Frankl-1978}, Frankl also made a conjecture on the maximum size of a $t$-intersecting family of ${[n]\choose k}$ for all positive integers $t,k$ and $n$. This conjecture was partially proved by Frankl and F\"{u}redi in  \cite{Frankl--Furedi-1991} and completely settled by Ahlswede and Khachatrian in \cite{Ahlswede-Khachatrian-1997}.

Determining the structure of extremal non-trivial $t$-intersecting families of ${[n]\choose k}$ was a long-standing problem. The first result was the Hilton-Milner Theorem \cite{Hilton-Milner-1967} which describes the structure of such families for $t=1$. A significant step was taken in \cite{Frankl-1978-1} by Frankl, who determined such families for $t\geq 2$ and $n>n_1(k,t)$. In \cite{Frankl-Furedi-1986}, Frankl and F\"{u}redi gave a short and elegant proof for the Hilton-Milner Theorem by using the shifting technique, and also asked whether $n_1(k,t)<ckt$ holds.  Ahlswede and Khachatrian \cite{Ahlswede-Khachatrian-1996} answered this question and gave a complete result on non-trivial intersection problems for finite sets. Recently, other maximal non-trivial $t$-intersecting families with large size had been studied. See \cite{Cao-set,Han-Kohayakawa,Kostochka-Mubayi}.

Let $r\geq 2$ be an integer. A family $\mathcal{F}\subseteq {[n]\choose k}$ is called \emph{$r$-wise $t$-intersecting} if for all $F_1,\ldots,F_r\in\mathcal{F}$ one has $|F_1\cap\cdots\cap F_r|\geq t$. An $r$-wise $t$-intersecting family is called \emph{trivial} if all its members contain a common specified $t$-subset of $[n]$, and \emph{non-trivial} otherwise. An $r$-wise $t$-intersecting family $\mathcal{F}$ is called \emph{maximal} if $\mathcal{F}\cup\{A\}$ is not $r$-wise $t$-intersecting for every $A\in{[n]\choose k}\setminus\mathcal{F}$. Observe that the `$2$-wise $t$-intersecting' is the classical `$t$-intersecting'.

From the Erd\H{o}s-Ko-Rado Theorem for finite sets, it is clear that for $n>(t+1)(k-t+1)$ each extremal $r$-wise $t$-intersecting families of ${[n]\choose k}$ is trivial. But for $n\leq(t+1)(k-t+1)$, it seems to be an extremely challenging task to completely determine the structure of extremal $r$-wise $t$-intersecting families. We refer the readers to
\cite{a1,Frankl-1976,Frankl-Tokushige-2002,Frankl-Tokushige-2005,Liu-Zhang-Li-Zhang-2016,Moura-1999,Tokushige-2006,Tokushige-2007,Tokushige-2007-2} for some known results. In \cite{O-V-2021}, O'Neill and Verstra\"{e}te characterized the structure of non-trivial $r$-wise $1$-intersecting families of ${[n]\choose k}$ with maximum size for large $n$, and gave a stability theorem about these extremal families. From their former result, the conjecture given by Hilton and Milner in \cite{Hilton-Milner-1967} was proved. In this paper, we will determine the structure of extremal non-trivial $r$-wise $t$-intersecting families for any $t$ and large $n$, and also gave a stability theorem. To state our theorems, we first introduce the following constructions of families.

Let $n$, $k$ and $d$ be positive integers with $n\geq 2k$ and $k>d$. Assume that $Z\in{[n]\choose d+2}$, $M\in{[n]\choose k+1}$ and $X\subseteq {M\choose d}$. Write
{\small\begin{align}
\mathcal{H}(k,d+1,X,M)=&\left\{H\in{[n]\choose k}\mid X\subseteq H,\ |H\cap M|\geq d+1\right\}\cup{M\choose k},\label{s-equ-2}\\
\mathcal{A}(k,d+1,Z)=&\left\{A\in{[n]\choose k}\mid |A\cap Z|\geq d+1\right\}.\label{s-equ-1}
\end{align}}%
For $c\in\{k+1,k+2,\ldots,2k-d,n\},$ set
{\small\begin{align}
h_1(d,k,c)=&{n-d\choose k-d}-{n-k\choose k-d}+{n-c\choose 2k-c-d}+d(c-k),\label{s-equ-4}\\
h_2(d+2)=&(d+2){n-d-1\choose k-d-1}-(d+1){n-d-2\choose k-d-2}.\label{s-equ-5}
\end{align}}%
It is routine to check that $|\mathcal{H}(k,d+1,X,M)|=h_1(d,k,k+1)$ and $|\mathcal{A}(k,d+1,Z)|=h_2(d+2)$. In Observation~\ref{s-obs1} (ii) and (iii) of Section 2, we will see that $h_1(d,k,n)$ is the size of the family
{\small
\begin{align*}
\left\{F\in{[n]\choose k}\mid X\subseteq F,\ |F\cap M|\geq d+1\right\}\cup\left\{F\in{[n]\choose k}\mid |F\cap X|=d-1,\ |F\cap M|=k-1\right\}.
\end{align*}
}%
Our first main theorem is as follows.
\begin{thm}\label{s-main-1-1}
	Let $n$, $k$, $t$ and $r$ be integers with $3\leq r$ and $(t+r-1)(k-t-r+3)<n$, and $\mathcal{F}\subseteq{[n]\choose k}$ a non-trivial $r$-wise $t$-intersecting family.
	\begin{enumerate}[{\rm(i)}]
\item If $2t+2r-3<k,$ then
		$$
		|\mathcal{F}|\leq\max\{h_1(t+r-2,k,k+1), h_2(t+r)\}
		$$
		and equality holds only if $\mathcal{F}$ is the one of $\mathcal{H}(k,t+r-1,X,M)$ and $\mathcal{A}(k,t+r-1,Z)$ with larger size, where $M$ is a $(k+1)$-subset of $[n]$, $X$ is an $(t+r-2)$-subset of $M$, and $Z$ is an $(r+t)$-subset of $[n]$.		
\item If $2t+2r-3\geq k$, then
		$$
		|\mathcal{F}|\leq h_2(t+r),
		$$
		and equality holds only if $\mathcal{F}$ is $\mathcal{A}(k,t+r-1,Z)$, where $Z$ is an $(t+r)$-subset of $[n]$.
	\end{enumerate}
\end{thm}

We mention here that it is difficult to completely compare the sizes of $h_2(t+r)$ and  $h_1(t+r-2,k,k+1)$ when $2t+2r-3<k.$ From \cite[Lemmas~2.6 and 5.2(ii)]{Cao-set}, we only know that $h_2(t+r)<h_1(t+r-2,k,k+1)$ if $2t+2r-1\leq k$ and $n\geq{t+r\choose 2}\cdot(k-t-r+3)^2+t+r-2$, or $2t+2r-2=k$ and $n$ is sufficiently large.

Next, we will give a stability result for Theorem~\ref{s-main-1-1}. When $k=t+r-1$, in Lemma~\ref{s-k=t+1} we will show that any non-trivial $r$-wise $t$-intersecting family of ${[n]\choose k}$ is a subfamily of ${M\choose k}$ for some $(k+1)$-subset $M$. In the following, we consider the case $k>t+r-1$. Assume that $n\geq{t+r\choose 2}\cdot(k-t-r+3)^2+t+r-2$ and $r\geq 3$. Then $t+r-2\geq t+1\geq 2$. By \cite[Remark~1, Lemmas~5.1 ans 5.2]{Cao-set}, observe that $h_1(t+r-2,k,k+1)>h_1(t+r-2,k,n)>h_2(t+r)$ if $4\leq 2t+2r-2<k$, $\min\{h_2(t+r), h_1(t+r-2,k,k+1)\}>h_1(t+r-2,k,n)$ if $k\leq 2t+2r-2\leq2k-4$, and $h_2(t+r)=h_1(t+r-2,k,n)>h_1(t+r-2,k,k+1)$ if $2t+2r-2=2k-2$.

\begin{thm}\label{s-main-1-2}
Let $r\geq 3$, $t+r\leq k$, and $\max\left\{{t+r\choose 2},\frac{k-t-r+4}{2}\right\}\cdot(k-t-r+3)^2+t+r-2\leq n$. Suppose that $\mathcal{F}\subseteq{[n]\choose k}$ is a non-trivial $r$-wise $t$-intersecting family.
\begin{enumerate}[{\rm(i)}]
\item If $2t+2r-2<k,$ and $|\mathcal{F}|>h_1(t+r-2,k,n)$, then $\mathcal{F}$ is a subfamily of $\mathcal{H}(k,t+r-1,X,M)$, where $M$ is a $(k+1)$-subset of $[n]$ and $X$ is an $(t+r-2)$-subset of $M.$

\item If $k\leq 2t+2r-2<2k-2$ and $|\mathcal{F}|>h_1(t+r-2,k,n)$, or $2t+2r-2=2k-2$ and $|\mathcal{F}|\geq h_1(t+r-2,k,k+1)$, then $\mathcal{F}$ is a subfamily of $\mathcal{H}(k,t+r-1,X,M)$ or $\mathcal{A}(k,t+r-1,Z)$, where $M$ is a $(k+1)$-subset of $[n]$, $X$ is an $(t+r-2)$-subset of $M$, and $Z$ is an $(t+r)$-subset of $[n]$.
\end{enumerate}
\end{thm}

Let $n$, $r$, $k_1,\ldots,k_r$ and $t$ be positive integers. We say that $r$ families $\mathcal{F}_1\subseteq{[n]\choose k_1},\ \mathcal{F}_2\subseteq{[n]\choose k_2},\ldots,\mathcal{F}_r\subseteq{[n]\choose k_r}$ are $r$-cross $t$-intersecting if $|F_1\cap F_2\cap\cdots\cap F_r|\geq t$	
for all $F_i\in\mathcal{F}_i,\ 1\leq i\leq r.$ The $r$-cross $t$-intersecting families $\mathcal{F}_1,\mathcal{F}_2,\ldots,\mathcal{F}_r$ are said to be \emph{trivial} if $|\cap_{1\leq i\leq r}\cap_{F\in\mathcal{F}_i}F|\geq t$ and \emph{non-trivial} otherwise, and said to be maximal if $\mathcal{G}_i=\mathcal{F}_i\ (i=1,2,\ldots,r)$ for any $r$-cross $t$-intersecting families $\mathcal{G}_1\subseteq{[n]\choose k_1},\ \mathcal{G}_2\subseteq{[n]\choose k_2},\ldots,\mathcal{G}_r\subseteq{[n]\choose k_r}$ with $\mathcal{F}_i\subseteq \mathcal{G}_i\ (i=1,2,\ldots,r)$. Observe that if $\mathcal{F}\subseteq {[n]\choose k}$ is an $r$-wise $t$-intersecting family, by setting $\mathcal{F}_1=\mathcal{F}_2=\cdots=\mathcal{F}_r=\mathcal{F}$, then $\mathcal{F}_1,\mathcal{F}_2,\ldots,\mathcal{F}_r$ are $r$-cross $t$-intersecting families. For convenience, we just say `cross $t$-intersecting' instead of `$2$-cross $t$-intersecting'.

There is a vast literature on studying the maximum sum or the maximum product of sizes of $r$-cross $t$-intersecting families for finite sets, and on studying their extremal structures. We refer the readers to \cite{Borg-2009, Frankl--Tokushige-1992, Wang-Zhang-2011} for the sum version, and \cite{Bey-2005,Borg-2014,Borg-2016,Frankl--Lee--Siggers--Tokushige-2014,Frankl--Tokushige-2011,Hilton-1977,Matsumoto-Tokushige-1989,Pyber-1986,Tokushige-2010,Tokushige-2013,Wang-Zhang-2011} for the product version. From these results, it can be deduced that for large $n$, if $\mathcal{F}_1,\ \mathcal{F}_2,\ldots,\mathcal{F}_r$ are $r$-cross $t$-intersecting families, then $\prod_{i=1}^r|\mathcal{F}_i|\leq\prod_{i=1}^r{n-t\choose k_i-t}$, and the extremal structure of these families is trivial. In this paper, we focus on the non-trivial cross $t$-intersecting families with maximum product of their sizes.

Let $\mathcal{F}_1\subseteq{[n]\choose k_1}$ and $\mathcal{F}_2\subseteq{[n]\choose k_2}$ be non-trivial cross $t$-intersecting families. Without loss of generality, assume that $k_1\geq k_2$. Suppose that $k_2=t$. Since $\mathcal{F}_1$ and $\mathcal{F}_2$ are non-trivial cross $t$-intersecting, we have $k_1\geq t+1,$ $|\mathcal{F}_2|\geq 2$ and $F_1\supseteq F_2$ for each $F_1\in\mathcal{F}_1$ and $F_2\in\mathcal{F}_2$. Let $W=\cup_{F_2\in\mathcal{F}_2}F_2$. Then $t+1\leq|W|\leq k_1,$ $\mathcal{F}_1\subseteq\{F_1\in{[n]\choose k_1}\mid W\subseteq F_1\}$ and $\mathcal{F}_2\subseteq{W\choose k_2}.$ Note that $\{F_1\in{[n]\choose k_1}\mid W\subseteq F_1\}$ and ${W\choose k_2}$ are cross $t$-intersecting. Therefore, for any non-trivial cross $t$-intersecting families $\mathcal{G}_1\subseteq{[n]\choose k_1}$ and $\mathcal{G}_2\subseteq{[n]\choose k_2}$ with $k_1\geq k_2=t$, there exists a subset $W$ such that $\mathcal{G}_1\subseteq\{G_1\in{[n]\choose k_1}\mid W\subseteq G_1\}$ and $\mathcal{G}_2\subseteq{W\choose k_2}$. The following theorem characterize the structure of $\mathcal{F}_1$ and $\mathcal{F}_2$ with maximum product of their sizes for $k_2\geq t+1$.

\begin{thm}\label{s-main-1}
Let $n$, $k_1$, $k_2$ and $t$ be positive integers with $n\geq\max\{t+1,\ k_2-t\}\cdot(t+1)(k_1-t+1)(k_2-t+1)+t+1$, $k_1\geq k_2\geq t+1$ and $(k_1,k_2,t)\neq (2,2,1),(3,2,1),(4,2,1)$ or $(4,4,2)$. Assume that $\mathcal{F}_1\subseteq{[n]\choose k_1}$ and $\mathcal{F}_2\subseteq{[n]\choose k_2}$ are a pair of non-trivial cross $t$-intersecting families.
\begin{enumerate}[{\rm(i)}]
\item If $k_2\geq 2t+1$, then
$$
|\mathcal{F}_1||\mathcal{F}_2|\leq \left({n-t\choose k_1-t}-{n-k_2-1\choose k_1-t}\right)\left({n-t\choose k_2-t}+t\right),
$$
and the equality holds only if there exist a $t$-subset $X$ and a $(k_2+1)$-subset $M$ of $[n]$ with $X\subseteq M$ such that one of the following holds:
\begin{itemize}
\item[{\rm(ia)}] $\mathcal{F}_1=\{F\in{[n]\choose k_1}\mid X\subseteq F,\ |F\cap M|\geq t+1\}$ and $\mathcal{F}_2=\{F\in{[n]\choose k_2}\mid X\subseteq F\}\cup {M\choose k_2}$,
\item[{\rm(ib)}] $\mathcal{F}_1=\{F\in{[n]\choose k_1}\mid X\subseteq F\}\cup {M\choose k_1}$ and $\mathcal{F}_2=\{F\in{[n]\choose k_2}\mid X\subseteq F,\ |F\cap M|\geq t+1\}$ with $k_1=k_2$.
\end{itemize}

\item If $t+1\leq k_2\leq 2t$, then
$$
|\mathcal{F}_1||\mathcal{F}_2|\leq {n-t-1\choose k_1-t-1}\left((t+1){n-t-1\choose k_2-t}+{n-t-1\choose k_2-t-1}\right),
$$
and the equality holds only if there exists a $(t+1)$-subset $T$ of $[n]$ with $X\subseteq M$ such that one of the following holds:
\begin{itemize}
\item[{\rm(ia)}] $\mathcal{F}_1=\{F\in{[n]\choose k_1}\mid T\subseteq F\}$ and $\mathcal{F}_2=\{F\in{[n]\choose k_2}\mid |F\cap T|\geq t\}$,
\item[{\rm(ib)}] $\mathcal{F}_1=\{F\in{[n]\choose k_1}\mid |F\cap T|\geq t\}$ and $\mathcal{F}_2=\{F\in{[n]\choose k_2}\mid T\subseteq F\}$ with $k_1=k_2$.
\end{itemize}
\end{enumerate}
\end{thm}

The rest of this paper is organized as follows. In Sections 2, we will prove Theorems~\ref{s-main-1-1} and \ref{s-main-1-2}. In Section 3, we will prove Theorems~\ref{s-main-1}. In Section 4, we will prove some inequalities used in this paper.

\section{Non-trivial $r$-wise $t$-intersecting families}
In this section, we describe the structure of non-trivial $r$-wise $t$-intersecting families with maximum size, and give a stability result for these families.

\subsection{Some properties}
\begin{lemma}\label{s-lem-1}
Let $\mathcal{F}\subseteq {[n]\choose k}$ be a non-trivial $r$-wise $t$-intersecting family. Suppose that $S$ is a subset of $[n]$ such that $|S\cap F_1\cap F_2\cap\cdots \cap F_{r-1}|\geq t$ for all $F_1,F_2,\ldots,F_{r-1}\in\mathcal{F}$. Then the following hold.
\begin{enumerate}[\rm(i)]
\item We have $|S|\geq t+r-1.$
\item Suppose that $A_1,A_2,\ldots,A_m$ are elements in $\mathcal{F}$ with $m\leq r-1$. Then $|S\cap A_1\cap A_2\cap \cdots\cap A_m|\geq r+t-1-m.$
\end{enumerate}
\end{lemma}
\proof (i)\quad Set $C=\cap_{F\in\mathcal{F}}F$, and then $|C|<t$ due to that $\mathcal{F}$ is non-trivial. Since $|S|\geq t$, we have $S\setminus C\neq\emptyset.$ Assume that $S\setminus C=\{x_1,x_2,\ldots,x_\ell\}$. Then there exist $A_1,A_2,\ldots,A_\ell\in\mathcal{F}$ such that $x_i\notin A_i$ for $1\leq i\leq \ell.$ It follows that $S\cap A_1\cap\cdots\cap A_{\ell}\subseteq C.$ If $\ell\leq r-1,$ then  $t>|C|\geq|S\cap A_1\cap A_2\cap\cdots\cap A_\ell|\geq t,$ a contradiction. Hence, we have $\ell> r-1.$ Note that $S\setminus\{x_1,\ldots,x_{r-1}\}\supseteq S\cap A_1\cap A_2\cap\cdots\cap A_{r-1}$, and $|S\cap A_1\cap A_2\cap\cdots\cap A_{r-1}|\geq t$, implying that $|S|-(r-1)\geq t.$ Therefore (i) holds.

(ii)\quad Observe that it is clear if $m=r-1$. Suppose that $m\leq r-2$ in the following. Since $|(S\cap A_1\cap A_2\cap\cdots\cap A_m)\cap F_1\cap\cdots\cap F_{r-1-m}|\geq t$ for all $F_1,F_2,\cdots,F_{r-1-m}\in\mathcal{F}$, and $\mathcal{F}$ is an $(r-m)$-wise $t$-intersecting family, by (i), we have $|S\cap A_1\cap\cdots\cap A_m|\geq t+r-1-m$ as required.    $\qed$

\begin{cor}\label{s-lem-2}
	If $r>k-t+1,$ then there does not exist a non-trivial $r$-wise $t$-intersecting  family of $k$-subsets of $[n]$.
\end{cor}
\proof Assume that there exists a non-trivial $r$-wise $t$-intersecting family $\mathcal{F}$ of $k$-subsets of $[n]$. Set $S$ be a $k$-subset in $\mathcal{F}$. Then $|S\cap F_1\cap F_2\cap\cdots \cap F_{r-1}|\geq t$ for all $F_1,F_2,\ldots,F_{r-1}\in\mathcal{F}$. By Lemma~\ref{s-lem-1}, we have $|S|=k\geq t+r-1$, a contradiction. $\qed$
\begin{cor}\label{s-lem-2-1-1}
Suppose that $B_1, B_2,\cdots,B_d$ are elements in the non-trivial $r$-wise $t$-intersecting family $\mathcal{F}\subseteq{[n]\choose k}$ with $d\leq r$.
Then $|B_1\cap B_2\cap\cdots\cap B_d|\geq t+r-d$. Moreover, $\mathcal{F}$ is a $(t+r-2)$-intersecting family.
\end{cor}
\proof Set $S=B_1$  in Lemma~\ref{s-lem-1} (ii), and then the former part of this corollary holds. In particular, set $d=2$, and then we have $\mathcal{F}$ is $(t+r-2)$-intersecting.   $\qed$

\subsection{The proofs of Theorems~\ref{s-main-1-1} and \ref{s-main-1-2}}

\begin{lemma} \label{s-family}
Let $M\in{[n]\choose k+1}$, $X\in{M\choose t+r-2}$ and  $Z\in{[n]\choose t+r}$. Then both $\mathcal{A}(k,t+r-1,Z)$ and $\mathcal{H}(k,t+r-1,X,M)$  are non-trivial $r$-wise $t$-intersecting families.
\end{lemma}
\proof Firstly, we claim that for positive integers $\ell$ and $m$ with $m\leq \ell,$ if $E_1,\ldots,E_m\in{E\choose\ell}$  for an $(\ell+1)$-set $E$, then $|E_1\cap E_2\cap\cdots\cap E_m|\geq\ell-m+1$. Indeed, it is clear for $m=1$; and for $m\geq 2,$
\begin{align*}
	|E_1\cap E_2\cap\cdots\cap E_m|\geq&|E_1\cap E_2\cap\cdots\cap E_{m-1}|+|E_m|-|E|\\
	=&|E_1\cap E_2\cap\cdots\cap E_{m-1}|-1\\
	\geq&|E_1\cap E_2\cap\cdots\cap E_{m-2}|-2\\
	&\vdots\\
	\geq& |E_1|-(m-1)\\
	=&\ell-m+1.
\end{align*}

Let $F_1,F_2,\ldots,F_r\in\mathcal{A}(k,t+r-1,Z).$ If $Z\subseteq F_i$ for every $i\in\{1,2,\ldots,r\}$, then $|F_1\cap\cdots\cap F_r|\geq t$. Now suppose that there exists $h\in\{1,2,\ldots,r\}$ such that $Z\nsubseteq F_i$. Without loss of generality, we can assume that
$$
|F_i\cap Z|=\left\{\begin{array}{ll} t+r-1, & i\leq h,\\t+r,& \mbox{otherwise.}\end{array}\right.
$$
Then, by the claim above, we have
\begin{align*}
	|F_1\cap\cdots\cap F_r|\geq&|(F_1\cap Z)\cap\cdots\cap (F_r\cap Z)|\\
	=&|(F_1\cap Z)\cap\cdots\cap (F_h\cap Z)|\\
	\geq&t+r-1-h+1\\
	\geq&t.
\end{align*}

Let $F_1,F_2,\ldots,F_r\in\mathcal{H}(k,t+r-1,X,M).$ If $F_i\in{M\choose k}$ for every $i\in\{1,2,\ldots,r\}$, then $|F_1\cap\cdots\cap F_r|\geq t$ due to $t+r-2\leq k+1$ and the claim above. Now suppose that there exists $h\in\{1,2,\ldots,r\}$ such that $F_i\notin{M\choose k}$. Without loss of generality, we can assume that
$$
F_i\in\left\{\begin{array}{ll} \mathcal{H}(k,t+r-1,X,M)\setminus{M\choose k}, & i\leq h,\\ {M\choose k},& \mbox{otherwise.}\end{array}\right.
$$
Then, if $h=1$, by the claim above and the construction of $\mathcal{H}(k,t+r-1,X,M)$, we have
\begin{align*}
	|F_1\cap\cdots\cap F_r|\geq& |((F_1\cap M)\cap F_2\cap\cdots\cap F_{r})|\\
	\geq&|F_1\cap M|+|F_2\cap\cdots\cap F_{r}|-|M|\\
	\geq&(t+r-1)+(k-r+2)-(k+1)\\
	\geq&t;
\end{align*}
if $h\geq 2$, by the claim above and the construction of $\mathcal{H}(k,t+r-1,X,M)$ again, we have
\begin{align*}
	|F_1\cap\cdots\cap F_r|\geq& |(F_1\cap\cdots\cap F_h\cap X)\cap (F_{h+1}\cap\cdots\cap F_r)|\\
	\geq&|F_1\cap\cdots\cap F_h\cap X|+|F_{h+1}\cap\cdots\cap F_r|-|M|\\
	\geq&(t+r-2)+(k-r+h+1)-(k+1)\\
	\geq&t.
\end{align*}

Hence, both $\mathcal{A}(k,t+r-1,Z)$ and $\mathcal{H}(k,t+r-1,X,M)$ are $r$-wise $t$-intersecting. It is clear that these two families are non-trivial, and therefore the lemma holds.  $\qed$

The following theorem is essential to prove Theorem~\ref{s-main-1-1}.
\begin{thm}{\rm (\cite{Ahlswede-Khachatrian-1996})}\label{s-non-trivial-1}
	Let $n$, $k$ and $r$ be integers with $(t+1)(k-t+1)<n$, and $\mathcal{F}\subseteq{[n]\choose k}$ a non-trivial $t$-intersecting family.
	\begin{enumerate}[{\rm(i)}]
		\item If $k\leq 2t+1$, then $|\mathcal{F}|\leq h_2(t+2),$
		and equality holds only if $\mathcal{F}=\mathcal{A}(k,t+1,Z)$ for some $(t+2)$-subset $Z$ of $[n]$.
		\item If $k> 2t+1,$ then $|\mathcal{F}|\leq\max\{h_1(t,k,k+1),\ h_2(t+2)\}$,
		and equality holds only if $\mathcal{F}$ is the one of $\mathcal{H}(k,t+1,X,M)$ and $\mathcal{A}(k,t+1,Z)$ with large size, where $M$ is a $(k+1)$-subset of $[n]$, $X$ is a $t$-subset of $M$, and $Z$ is a $(t+2)$-subset of $[n]$.
	\end{enumerate}
\end{thm}

\noindent\textbf{The proof of Theorem~\ref{s-main-1-1}.}\quad Since $\mathcal{F}$ is non-trivial $r$-wise $t$-intersecting family, by Corollary~\ref{s-lem-2-1-1}, $\mathcal{F}$ is a non-trivial $(t+r-2)$-intersecting family. By Theorem~\ref{s-non-trivial-1} and Lemma~\ref{s-family}, the required result holds. $\qed$
\begin{lemma}\label{s-k=t+1}
Let $n$, $k$, $t$ and $r$ be positive integers with $r\geq 2$, $k=t+r-1$ and $n\geq k+1.$ If $\mathcal{F}\subseteq{[n]\choose k}$ is a maximal non-trivial $r$-wise $t$-intersecting family, then $\mathcal{F}={M\choose k}$ for some $(k+1)$-subset $M$ of $[n]$.
\end{lemma}
\proof By Corollary~\ref{s-lem-2-1-1}, we have that $\mathcal{F}$ is a non-trivial $(k-1)$-intersecting family. It is well known from the structure of the cliques of Johnson graphs that $\mathcal{F}\subseteq{M\choose k}$ for some $(k+1)$-subset $M$ of $[n]$. From the claim in the proof of Lemma~\ref{s-family}, observe that ${M\choose k}$ is $r$-wise $t$-intersecting. It follows from the maximality of $\mathcal{F}$ that $\mathcal{F}={M\choose k}$ holds.    $\qed$

To prove Theorem~\ref{s-main-1-2}, let us first introduce the following family.

\noindent\textbf{Family I.}\quad Let $X$, $M$ and $C$ be three subsets of $[n]$ such that $X\subseteq M\subseteq C$, $|X|=t$, $|M|=k$ and $|C|=c$, where $k\geq t+1$, $n\geq 2k$ and $c\in\{k+1,k+2,\ldots,2k-t,n\}$. Denote
{\small\begin{align}\label{s-eq:H2}
	 \mathcal{H}_1(X,M,C)=\mathcal{E}_1(X,M)\cup\mathcal{E}_2(X,M,C)\cup\mathcal{E}_2(X,M,C),
\end{align}}%
where
{\small\begin{align*}
	\mathcal{E}_1(X,M)=&\left\{F\in{[n]\choose k}\mid X\subseteq F,\ |F\cap M|\geq t+1\right\},\\
	\mathcal{E}_2(X,M,C)=&\left\{F\in{[n]\choose k}\mid F\cap M=X,\ |F\cap C|=c-k+t\right\},\\
	\mathcal{E}_3(X,M,C)=&\left\{F\in{C\choose k}\mid |F\cap X|=t-1,\ |F\cap M|=k-1\right\}.
\end{align*}}%

By this construction and \cite[Remark~1, Lemma~2.5]{Cao-set}, the following observation holds.

\begin{obs}\label{s-obs1} {\rm Suppose $X$, $M$ and $C$ are three subsets of $[n]$ satisfying the condition in Family I. Then the following hold.
\begin{itemize}
\item[{\rm(i)}] If $|C|=k+1$, then $\mathcal{H}_1(X,M,C)=\mathcal{H}(k,t+1,X,C).$		
\item[{\rm(ii)}] $|\mathcal{H}_1(X,M,C)|=h_1(t,k,c)$, where $c=|C|\in\{k+1,k+2,\ldots,2k-t,n\}$.
\item[{\rm(iii)}] $\mathcal{H}_1(X,M,[n])=\mathcal{E}_1(X,M)\cup \mathcal{E}_3(X,M,[n])$. If $t=k-2$, then $\mathcal{H}_1(X,M,[n])=\mathcal{A}(k,k-1,M)$ and $h_1(k-2,k,n)=h_2(k)$.
\item[{\rm(iv)}] For each $F\in\mathcal{E}_2(X,M,C)$, since $|C\setminus M|=c-k$, and $|F\cap C|=|X|+|F\cap(C\setminus M)|=c-k+t$, we have $|F\cap(C\setminus M)|=c-k$ and $C\setminus M\subseteq F$.
\end{itemize} }
\end{obs}

The following theorem given in \cite{Cao-set} characterized the structure of maximal non-trivial uniform $t$-intersecting  families with large size for finite sets.
\begin{thm}{\rm(\cite[Theorem~1.1]{Cao-set})}\label{s-main-1-non-trivial}
	Let $1\leq t\leq k-2$, and $\max\left\{{t+2\choose 2},\frac{k-t+2}{2}\right\}\cdot(k-t+1)^2+t\leq n$. If  $\mathcal{F}\subseteq{[n]\choose k}$ is a maximal non-trivial $t$-intersecting family with
	$$
	|\mathcal{F}|\geq (k-t){n-t-1\choose k-t-1}-{k-t\choose 2}{n-t-2\choose k-t-2},
	$$
	then one of the following holds.
	\begin{itemize}
		\item[{\rm(i)}] $\mathcal{F}=\mathcal{H}_1(X,M,C)$ for some $t$-subset $X$, $k$-subset $M$ and $c$-subset $C$ of $[n]$ with $X\subseteq M\subseteq C$, where $c\in\{k+1,k+2,\ldots,2k-t,n\}$.
		\item[{\rm(ii)}] $\mathcal{F}=\mathcal{A}(k,t+1,Z)$ for some $(t+2)$-subset $Z$ of $[n],$ and $\frac{k}{2}-1\leq t\leq k-2$.
	\end{itemize}
\end{thm}
\begin{lemma}\label{s-wise-sub}
Let $n$, $k$, $t$ and $r$ be positive integers with $r\geq 3$, $t+r-2\leq k-2$ and $2k+t+r+2\leq n$. Let $\mathcal{F}$ be a non-trivial $r$-wise $t$-intersecting subfamily of $\mathcal{H}_1(X,M,C)$, where $X\subseteq M\subseteq C,$ $|X|=t+r-2$, $|M|=k$ and $|C|=c\in\{k+2,k+3,\ldots, 2k-t-r+2\}$. If $|\mathcal{F}|> h_1(t+r-2,k,n),$ then $\mathcal{F}\subseteq\mathcal{H}(k,t+r-1,X,M^\prime)$ for some $(k+1)$-subset $M^{\prime}$ of $C$ with $M\subseteq M^{\prime}$.
\end{lemma}
\proof  Let $\mathcal{E}_1(X,M),$ $\mathcal{E}_2(X, M,C)$ and $\mathcal{E}_3(X, M, C)$ be as in (\ref{s-eq:H2}). Observe that all the subsets in $\mathcal{E}_1(X,M)$ and $\mathcal{E}_2(X, M,C)$ are containing $X$. Since $\mathcal{F}$ is non-trivial, we have $\mathcal{F}\cap \mathcal{E}_3(X, M, C)\neq\emptyset.$

Since $\mathcal{F}\subseteq \mathcal{H}_1(X,M,C)$ and $|\mathcal{F}|> h_1(t+r-2,k,n)$, by Observation~\ref{s-obs1} (iii),
we have
\begin{align*}
	|\mathcal{F}\cap\mathcal{E}_2(X, M,C)|\geq& |\mathcal{F}|-|\mathcal{E}_1(X,M)|-|\mathcal{E}_3(X, M, C)|\\
	>& |\mathcal{H}_1(X,M,[n])|-|\mathcal{E}_1(X,M)|-|\mathcal{E}_3(X, M, C)|\\
	=&|\mathcal{E}_3(X, M, [n])|-|\mathcal{E}_3(X, M, C)|>0,
\end{align*}
which implies that $\mathcal{F}\cap\mathcal{E}_2(X, M,C)\neq\emptyset$. Choose $F^\prime\in\mathcal{F}\cap\mathcal{E}_2(X, M,C).$

Suppose for the contrary that $\mathcal{F}\nsubseteq \mathcal{H}(k,t+r-1,X,M^\prime)$ for each $(k+1)$-subset $M^\prime$ of $C$ with $M\subseteq M^\prime.$

\medskip

\textbf{Claim~1.} There exist $F_1,\;F_2\in\mathcal{F}\cap \mathcal{E}_3(X, M, C)$ satisfying $F_1\cap X\neq F_2\cap X$ and $F_2\nsubseteq F_1\cup M.$

Choose $G_1\in \mathcal{F}\cap \mathcal{E}_3(X, M, C)$. We have $G_1\cup M\subseteq C,$ and $|G_1\cup M|=k+1$ from $|G_1\cap M|=k-1.$  Consider the construction of $\mathcal{E}_1(X,M)$, $\mathcal{E}_2(X,M,C)$ and $\mathcal{H}(k,t+r-1, X, G_1\cup M)$. Observe that $\mathcal{E}_1(X,M)\subseteq\mathcal{H}(k,t+r-1, X, G_1\cup M)$. For each $F\in\mathcal{E}_2(X,M,C)$, since $(G_1\cup M)\setminus M\subseteq C\setminus M\subseteq F$ from Observation~\ref{s-obs1} (iv), we have $|F\cap (G_1\cup M)|\geq t+r-1$ due to $F\cap M=X$ and $|(G_1\cup M)\setminus M|\geq 1$, implying that $F\in\mathcal{H}(k,t+r-1, X, G_1\cup M)$. Hence $\mathcal{E}_2(X,M,C)\subseteq \mathcal{H}(k,t+r-1, X, G_1\cup M)$. If $\mathcal{F}\cap\mathcal{E}_3(X,M,C)\subseteq {G_1\cup M\choose k}$, then $\mathcal{F}\cap(\mathcal{E}_1(X,M)\cup \mathcal{E}_2(X,M,C)\cup \mathcal{E}_3(X,M,C))\subseteq\mathcal{H}(k,t+r-1, X, G_1\cup M),$ a contradiction. Hence, we have $\mathcal{F}\cap\mathcal{E}_3(X,M,C)\nsubseteq {G_1\cup M\choose k}$, implying that
there exists $G_2\in\mathcal{F}\cap\mathcal{E}_3(X,M,C)$ such that $G_2\nsubseteq G_1\cup M.$

If $G_2\cap X\neq G_1\cap X$, setting $F_1=G_1$ and $F_2=G_2$, then $F_1$ and $F_2$ are the required subsets.

Now assume that $G_3\cap X=G_1\cap X$ for each $G_3\in\mathcal{F}\cap\mathcal{E}_3(X,M,C)$ with $G_3\nsubseteq G_1\cup M$. Observe that $X\subseteq F$ for each $F\in\mathcal{F}\cap(\mathcal{E}_1(X,M)\cup\mathcal{E}_2(X,M,C))$, and $|F\cap X|=t+r-3\geq t$ for each $F\in\mathcal{F}\cap\mathcal{E}_3(X,M,C)$. Since $\mathcal{F}$ is non-trivial, there exists $G_4\in\mathcal{F}\cap\mathcal{E}_3(X,M,C)$ such that $G_4\cap X\neq G_1\cap X=G_2\cap X$. It follows from the assumption that $G_4\in{G_1\cup M\choose k}$, and so $G_4\cup M=G_1\cup M$ due to $|G_4\cap M|=k-1$ and $|G_1\cup M|=k+1$. Set $F_1=G_4$ and $F_2=G_2$. We have that $F_1$ and $F_2$ are the required subsets, and Claim 1 holds.

\medskip

By Claim~1, let $F_1$ and $F_2$ be the subsets in $\mathcal{F}\cap\mathcal{E}_3(X,M,C)$ with  $F_1\cap X\neq F_2\cap X$ and $F_2\nsubseteq F_1\cup M.$ If there exists $y\in F_1\cap F_2$ such that $y\notin M$, then $F_1\cup M=\{y\}\cup M=F_2\cup M,$ a contradiction. Thus $F_1\cap F_2\subseteq M\cap F_1\cap F_2$, implying that $F_1\cap F_2=M\cap F_1\cap F_2.$ Recall that $\mathcal{F}\cap\mathcal{E}_2(X, M,C)\neq\emptyset$. Then we can choose $F^\prime\in\mathcal{F}\cap\mathcal{E}_2(X, M,C).$ Observe that $F^\prime \cap F_1\cap F_2=F^\prime\cap M\cap F_1\cap F_2= X \cap F_1\cap F_2$ and $|X \cap F_1\cap F_2|=|X \cap F_1|+|X \cap F_2|-|(X \cap F_1)\cup (X \cap F_2)|=t+r-4$.

If $r=3$, then $|F^\prime \cap F_1\cap F_2|=|X \cap F_1\cap F_2|=t-1$, which is impossible due to that $\mathcal{F}$ is non-trivial $r$-wise $t$-intersecting. Now suppose that $r\geq 4$. For each subset $S$ of $X$ with $|S|\geq t$, since $\mathcal{F}$ is a non-trivial $r$-wise $t$-intersecting subfamily, there exists $F\in \mathcal{F}\cap \mathcal{E}_3(X, M, C)$ with $S\nsubseteq F,$ and furthermore, $|S\cap F|=|S|-1$ from $|F\cap X|=|X|-1$.
 Using this fact repeatedly, we can get $F_3,\ldots,F_{r-1}\in\mathcal{F}\cap\mathcal{E}_3(X,M,C)$ such that $|(X\cap F_1\cap F_2\cap\cdots\cap F_{i-1})\cap F_i|=|X\cap F_1\cap F_2\cap\cdots\cap F_{i-1}|-1$ for each $i\in\{3,\ldots,r-1\}$.
 Hence $|X\cap F_1\cap F_2\cap\cdots\cap F_{r-1}|=t-1,$ implying that $|F^\prime\cap F_1\cap F_2\cap\cdots\cap F_{r-1}|=t-1,$ a contradiction.    $\qed$

\noindent\textbf{The proof of Theorem~\ref{s-main-1-2}}\quad Since $\mathcal{F}$ is a non-trivial $r$-wise $t$-intersecting family, by Corollary~\ref{s-lem-2-1-1}, $\mathcal{F}$ is a non-trivial $(t+r-2)$-intersecting family. Suppose that $|\mathcal{F}|>h_1(t+r-2,k,n)$ if $t+r-2<k-2$, and $|\mathcal{F}|\geq h_1(t+r-2,k,k,k+1)$ if $t+r-2=k-2.$ Since ${t+r\choose 2}(k-t-r+3)^2+t+r-2\leq n$, by \cite[Lemma~2.6 (i)]{Cao-set}, we have that
\begin{align*}
|\mathcal{F}|>(k-t-r+2){n-t-r+1\choose k-t-r+1}-{k-t-r+2\choose 2}{n-t-r\choose k-t-r}.
\end{align*}
It follows from Theorem~\ref{s-main-1-non-trivial} that one of the following holds:
	\begin{enumerate}[{\rm(i)}]
	\item $\mathcal{F}\subseteq\mathcal{H}_1(X,M,C)$ for some $(t+r-2)$-subset $X$, $k$-subset $M$ and $c$-subset $C$ of $[n]$ with $X\subseteq M\subseteq C$, where $c\in\{k+1,k+2,\ldots,2k-t-r+2,n\}$,
	\item $\mathcal{F}\subseteq\mathcal\mathcal{A}(k,t+2-1,X,M)$ for some $(t+r)$-subset $Z$ of $[n],$ and $\frac{k}{2}-1\leq t+r-2\leq k-2$.
	\end{enumerate}

If $t+r-2<k-2$, by $|\mathcal{F}|>h_1(t+r-2,k,n)$, then $\mathcal{F}$ is not a subfamily of $\mathcal{H}_1(X,M,[n])$, implying that Theorem~\ref{s-main-1-2} holds due to Lemma~\ref{s-wise-sub}. Now suppose $t+r-2=k-2$. Since $|\mathcal{F}|\geq h_1(t+r-2,k,k,k+1)$, $\mathcal{F}$ is not a subfamily of $\mathcal{H}_1(X,M,C)$ with $c\in\{k+2,\ldots,2k-t-r+2\}$ due to \cite[Lemma~5.1 (i)]{Cao-set}, which showed that $h_1(t,k,c)$ decreases as $c\in\{k+1,\ldots,2k-t-r+2\}$ increases.  Thus Theorem~\ref{s-main-1-2} holds for $t+r-2=k-2$ from
Observation~\ref{s-obs1} (i) (iii).  $\qed$

\section{Non-trivial cross $t$-intersecting families}

In this section, we first study the structure or the product of the sizes of cross $t$-intersecting families with given $t$-covering numbers (Subsections~3.1 and 3.2), and then prove Theorem~\ref{s-main-1} (Subsection~3.3).

For any $\mathcal{F}\subseteq{[n]\choose k}$ (not necessary to be $t$-intersecting), a subset $T$ of $[n]$ is called to be a $t$-\emph{cover} of $\mathcal{F}$ if $|T\cap F|\geq t$ for all $F\in\mathcal{F}$, and the $t$-\emph{covering number} $\tau_t(\mathcal{F})$ of $\mathcal{F}\subseteq {[n]\choose k}$ is the minimum size of a $t$-cover of $\mathcal{F}$. Let $\mathcal{F}\subseteq{[n]\choose k}$ and ${\mathcal{G}\subseteq{[n]\choose \ell}}$ be cross $t$-intersecting families. It is clear that each element in $\mathcal{G}$ is a $t$-cover of $\mathcal{F}$, and each element in $\mathcal{F}$ is a $t$-cover of $\mathcal{G}$. Hence, $t\leq \tau_t(\mathcal{F})\leq \ell$ and $t\leq \tau_t(\mathcal{G})\leq k$.%
\begin{lemma}\label{s-tt}
Let $n$, $k$, $\ell$ and $t$ be positive integers with $n\geq k+\ell$, and $\mathcal{F}\subseteq{[n]\choose k}$ and $\mathcal{G}\subseteq{[n]\choose \ell}$ be maximal cross $t$-intersecting families. Assume that $\mathcal{T}_f$ is the set of the $t$-covers of $\mathcal{F}$ with size $\tau_t(\mathcal{F})$, and $\mathcal{T}_g$ is the set of the $t$-covers of $\mathcal{G}$ with size $\tau_t(\mathcal{G})$. Then $\mathcal{T}_f$ and $\mathcal{T}_g$ are cross $t$-intersecting families.
\end{lemma}
\proof Let $T_f\in \mathcal{T}_f$ and $T_g\in\mathcal{T}_g.$ Since $\tau_t(\mathcal{F})\leq \ell$ and $\tau_t(\mathcal{G})\leq k$, we have $|T_f\cup T_g|\leq \tau_t(\mathcal{F})+\tau_t(\mathcal{G})$. It follows from $n\geq k+\ell$ that $$
|[n]\setminus(T_f\cup T_g)|\geq n-\tau_t(\mathcal{F})-\tau_t(\mathcal{G})\geq (\ell-\tau_t(\mathcal{F}))+(k-\tau_t(\mathcal{G})).
$$
Then there exist $\ell$-subset $G^\prime$ and $k$-subset $F^\prime$ of $[n]$ such that $T_f\subseteq G^\prime$, $T_g\subseteq F^\prime$ and $F^\prime\cap G^\prime=T_f\cap T_g.$
For each $F\in\mathcal{F}$, observe that $|F\cap G^\prime|\geq t$ due to $|F\cap T_f|\geq t$ and $T_f\subseteq G^\prime$; similarly, observe that $|F^\prime\cap G|\geq t$ for each $G\in\mathcal{G}$. By the maximality of $\mathcal{F}$ and $\mathcal{G},$ we have $F^\prime\in\mathcal{F}$ and $G^\prime\in\mathcal{G}$. Since $\mathcal{F}$ and $\mathcal{G}$ are cross $t$-intersecting families, and $F^\prime\cap G^\prime=T_f\cap T_g$, we have $|T_f\cap T_g|=|F^\prime\cap G^\prime|\geq t$. Hence, $\mathcal{T}_f$ and $\mathcal{T}_g$ are cross $t$-intersecting families.  $\qed$

Suppose that $\mathcal{F}\subseteq{[n]\choose k}$ and $\mathcal{G}\subseteq{[n]\choose \ell}$ are maximal cross $t$-intersecting families. By Lemma~\ref{s-tt}, we have that $\mathcal{F}$ and $\mathcal{G}$ are non-trivial if and only if $\tau_t(\mathcal{F})\geq t+1$ or $\tau_t(\mathcal{G})\geq t+1.$ In the following proofs, for the sake of simplicity, we write
{\small\begin{align}
g_1(k,\ell,n,t)=&\left({n-t\choose k-t}-{n-\ell-1\choose k-t}\right)\left({n-t\choose \ell-t}+t\right),\label{s-eq-1}\\
g_2(k,\ell,n,t)=&{n-t-1\choose k-t-1}\left((t+1){n-t-1\choose \ell-t}+{n-t-1\choose\ell-t-1}\right),\label{s-eq-2}\\
g_3(k,\ell,n,t)=&\left((\ell-t){n-t-1\choose k-t-1}+{n-t-2\choose k-t-2}\right)\left({n-t\choose \ell-t}+t(k-t){n-t-2\choose \ell-t-1}\right),\label{s-eq-3}\\
g_4(k,\ell,n,t)=&(k-t+1)^2{t+2\choose 2}{n-t\choose k-t}{n-t-2\choose \ell-t-2},\label{s-eq-4}\\
g_5(k,\ell,n,t)=&(t+1)^2(k-t+1)(\ell-t+1){n-t-1\choose k-t-1}{n-t-1\choose \ell-t-1}.\label{s-eq-5}
\end{align}}%
Note that the upper bounds in Theorem~\ref{s-main-1} (i) and (ii) are $g_1(k_1,k_2,n,t)$ and $g_2(k_1,k_2,n,t)$, respectively. For a family $\mathcal{F}\subseteq {[n]\choose k}$ and a subset $S$ of $[n]$, define $\mathcal{F}_S=\left\{F\in\mathcal{F}\mid S\subseteq F\right\}$.

\subsection{$\{\tau_t(\mathcal{F}), \tau_t(\mathcal{G})\}=\{t,t+1\}$}

\begin{lemma}\label{s-FGtt+1}
Let $n$, $k$, $\ell$ and $t$ be positive integers with $n\geq 2(k-t-1)(\ell+1-t)+t+1,$ and $\mathcal{F}\subseteq{[n]\choose k}$ and $\mathcal{G}\subseteq{[n]\choose \ell}$ be maximal cross $t$-intersecting families with $\tau_t(\mathcal{F})=t$ and $\tau_t(\mathcal{G})=t+1.$ Assume that $X$ is a $t$-cover of $\mathcal{F}$ with size $t$, and $T$ is a $t$-cover of $\mathcal{G}$ with size $t+1.$
\begin{enumerate}[{\rm(i)}]
\item For each $F\in \mathcal{F}$ and $G\in\mathcal{G}\setminus\mathcal{G}_X,$ we have $X\subseteq F$, $|G\cap X|=t-1$ and $|F\cap(G\cup X)|\geq t+1.$
\item If there exists $(\ell+1)$-subset $M$ of $[n]$ such that $G\cup X=M$ for each $G\in\mathcal{G}\setminus\mathcal{G}_X$, then
$$
\mathcal{F}=\left\{F\in{[n]\choose k}\mid X\subseteq F,\ |F\cap M|\geq t+1\right\}\ \mbox{and}\ \mathcal{G}=\left\{G\in{[n]\choose \ell}\mid X\subseteq G\right\}\cup {M\choose \ell}.
$$
Moreover, $|\mathcal{F}||\mathcal{G}|=g_1(k,\ell,n,t)$.
\item If $T\subseteq F$ for each $F\in\mathcal{F}$, then
$$
\mathcal{F}=\left\{F\in{[n]\choose k}\mid T\subseteq F\right\}\quad\mbox{and}\quad \mathcal{G}=\left\{G\in{[n]\choose \ell}\mid |G\cap T|\geq t\right\}.
$$
Moreover, $|\mathcal{F}||\mathcal{G}|=g_2(k,\ell,n,t)$.
\item If there exist $G_1, G_2\in \mathcal{G}\setminus\mathcal{G}_X$ such that $G_1\cup X\neq G_2\cup X,$ and there exists $F^\prime\in\mathcal{F}$ such that $T\nsubseteq F^\prime,$ then $|\mathcal{F}||\mathcal{G}|\leq g_3(k,\ell,n,t)$.
\end{enumerate}
\end{lemma}
\proof (i)\quad For each $F\in\mathcal{F}$, it is clear that $X\subseteq F$ since $X$ is a $t$-cover of $\mathcal{F}$ with size $t$. By Lemma~\ref{s-tt}, note that $X\subseteq T$, and so $|T\setminus X|=1$.  For each $G\in\mathcal{G}\setminus\mathcal{G}_X$, since $|G\cap T|\geq t$, $|T\setminus X|=1$ and $X\nsubseteq G$, we have
$$
t-1\leq|G\cap T|-1\leq|G\cap X|\leq t-1,
$$
implying that $|G\cap X|=t-1$. Assume that $X=(G\cap X)\cup\{y_0\}$. Then $G\cup X=G\cup\{y_0\}$, and so
$$
|F\cap(G\cup X)|=|F\cap G|+|F\cap\{y_0\}|\geq t+1.
$$

(ii)\quad Set
\begin{align*}
\mathcal{F}_1=& \left\{F\in{[n]\choose k}\mid X\subseteq F,\ |F\cap M|\geq t+1\right\},\\
\mathcal{G}_1=&\left\{G\in{[n]\choose \ell}\mid X\subseteq G\right\}\cup {M\choose \ell}.
\end{align*}
It is routine to check that $\mathcal{F}_1$ and $\mathcal{G}_1$ are cross $t$-intersecting families. For every $F\in \mathcal{F},$ by (i), we have $F\in\mathcal{F}_1$. For every $G\in\mathcal{G}$, if $X\subseteq G$, then $G\in\mathcal{G}_1;$ if $X\nsubseteq G,$ then $G\cup X=M$, implying that $G\in{M\choose \ell}\subseteq\mathcal{G}_1$. By the maximality, we have that $\mathcal{F}=\mathcal{F}_1$ and $\mathcal{G}=\mathcal{G}_1.$ Observe that
\begin{align*}
\mathcal{F}_1=& \left\{F\in{[n]\choose k}\mid X\subseteq F\right\}\setminus \left\{F\in{[n]\choose k}\mid F\cap M=X\right\},\\
\mathcal{G}_1=&\left\{G\in{[n]\choose \ell}\mid X\subseteq G\right\}\cup \left\{G\in{M\choose \ell}\mid X\nsubseteq G\right\}.
\end{align*}
Then it is routine to check that $|\mathcal{F}||\mathcal{G}|=g_2(k,\ell,n,t)$ holds.

(iii)\quad  Observe that $\mathcal{F}\subseteq \left\{F\in{[n]\choose k}\mid T\subseteq F\right\}$ and $\mathcal{G}\subseteq\left\{G\in{[n]\choose \ell}\mid |G\cap T|\geq t\right\}$. Since $\left\{F\in{[n]\choose k}\mid T\subseteq F\right\}$ and $\left\{G\in{[n]\choose \ell}\mid |G\cap T|\geq t\right\}$ are cross $t$-intersecting, we have that the former part of (iii) holds by the maximality. Furthermore, it is routine to check that $|\mathcal{F}||\mathcal{G}|=g_2(k,\ell,n,t)$ holds.

(iv)\quad Assume that $G_1\cup X=M_1$, $G_2\cup X=M_2$, $M^\prime=M_1\cap M_2$ and $|M^\prime|=m.$ Observe that $t\leq m\leq\ell.$ Since $X\subseteq F$ for each $F\in\mathcal{F}$, we have $F^\prime\cap T=X$. Set
\begin{align*}
	\mathcal{F}^\prime=&\left\{F\in{[n]\choose k}\mid X\subseteq F,\ |F\cap M_1|\geq t+1,\ |F\cap M_2|\geq t+1\right\}, \\
	\mathcal{G}^\prime=&\left\{G\in{[n]\choose \ell}\mid X\nsubseteq G,\ |G\cap T|=t,\ |G\cap F^\prime|\geq t\right\}.
\end{align*}
By (i) and the assumption, we have that $\mathcal{F}\subseteq\mathcal{F}^\prime$ and $\mathcal{G}\subseteq\{G\in{[n]\choose \ell}\mid X\subseteq G\}\cup \mathcal{G}^\prime.$ It is clear that $|\{G\in{[n]\choose \ell}\mid X\subseteq G\}|={n-t\choose \ell-t}$. In order to prove that (iv) holds, it suffices to give the upper bounds of $|\mathcal{F}^\prime|$ and $|\mathcal{G}^\prime|$.

\;

\noindent{\textbf{Step 1.}}\quad Show that $|\mathcal{F}^\prime|\leq (\ell-t){n-t-1\choose k-t-1}+{n-t-2\choose k-t-2}$.

For each $m\in\{t,t+1,\ldots,\ell\}$, write
$$
f^{\prime}(n,k,\ell,m,t)=(m-t){n-t-1\choose k-t-1}+(\ell+1-m)^2{n-t-2\choose k-t-2}.
$$
 Observe that
\begin{align*}
f^{\prime}(n,k,\ell,m,t)=&(m-t)\left({n-t-1\choose k-t-1}-2(\ell+1-t){n-t-2\choose k-t-2}\right)\\&+\left((\ell+1-t)^2+(m-t)^2\right){n-t-2\choose k-t-2}.
\end{align*}
Since $n\geq 2(k-t-1)(\ell+1-t)+t+1,$  we have
$$
{n-t-1\choose k-t-1}-2(\ell+1-t){n-t-2\choose k-t-2}=\left(\frac{n-t-1}{k-t-1}-2(\ell+1-t)\right){n-t-2\choose k-t-2}\geq0,
$$
implying that $f^{\prime}(n,k,\ell,m,t)$ increases as $m\in\{t,t+1,\ldots,\ell\}$ increases.

Set
$$
\mathcal{W}=\left\{(W_1,W_2)\in{M_1\choose t+1}\times{M_2\choose t+1}\mid X\subseteq W_1\nsubseteq M^\prime,\ X\subseteq W_2\nsubseteq M^\prime \right\}.
$$
If there exists $(W_1,W_2)\in\mathcal{W}$ such that $W_1=W_2,$ then $W_1=W_2\subseteq M_1\cap M_2=M^\prime$, a contradiction. Hence, $|W_1\cap W_2|=t$ and $|W_1\cup W_2|=t+2$ for each $(W_1,W_2)\in\mathcal{W}.$

For each $F\in\mathcal{F}^\prime,$ if $|F\cap M^\prime|\geq t+1,$ then there exists $H\in{M^\prime \choose t+1}$ such that $X\subseteq H\subseteq F;$ if $|F\cap M^\prime|=t,$ then there exist $(W_1,W_2)\in\mathcal{W}$ such that $W_1\subseteq F$ and $W_2\subseteq F$ due to $|F\cap M_1|\geq t+1$ and $|F\cap M_2|\geq t+1$. Therefore, we have
$$
\mathcal{F}^\prime\subseteq \left(\bigcup_{H\in{M^\prime\choose t+1},\ X\subseteq H}\mathcal{F}^\prime_H\right)\bigcup\left(\bigcup_{(W_1,W_2)\in\mathcal{W}}\mathcal{F}^\prime_{W_1\cup W_2}\right).
$$
Observe that $|\mathcal{F}^\prime_{H}|\leq{n-t-1\choose k-t-1}$ for each $H\in{M^\prime \choose t+1}$ with $X\subseteq H$, $|\mathcal{W}|=(\ell+1-m)^2$ and $|\mathcal{F}^\prime_{W_1\cup W_2}|\leq{n-t-2\choose k-t-2}$ for each $(W_1,W_2)\in\mathcal{W}$. Then $|\mathcal{F}^\prime|\leq f^{\prime}(n,k,\ell,m,t)\leq f^{\prime}(n,k,\ell,\ell,t)=(\ell-t){n-t-1\choose k-t-1}+{n-t-2\choose k-t-2}$ since $f^{\prime}(n,k,\ell,m,t)$ increases as $m\in\{t,t+1,\ldots,\ell\}$.

\;

\noindent{\textbf{Step 2.}}\quad Show that $|\mathcal{G}^\prime|\leq t(k-t){n-t-2\choose \ell-t-1}$.

Set
$$
\mathcal{W}^\prime=\left\{W\in{T\cup F^\prime\choose t+1}\mid |W\cap T|=t,\ W\cap T\neq X \right\}.
$$
Observe that $|T\cup F^\prime|=k+1$ from $T\cap F^\prime=X$, and $|\mathcal{W}^\prime|=t(k-t).$ Let $G\in\mathcal{G}^\prime$. From $X\nsubseteq G$ and $|G\cap T|=t$, there exist $y_1\in X$ such that $G\cap T=T\setminus\{y_1\}$. Since $|G\cap F^\prime|\geq t$ and $F^\prime \cap T=X$, there exist $z_1\in\mathcal{F}^\prime\setminus T$ such that $z_1\in G$. Set $w=(T\setminus\{y_1\})\cup\{z_1\}$, and then $w\in\mathcal{W}^\prime$ and $W\subseteq G.$ Therefore,
$$
\mathcal{G}^\prime\subseteq\bigcup_{W\in\mathcal{W}^\prime}\mathcal{G}^\prime_{W}\subseteq \bigcup_{W\in\mathcal{W}^\prime}\left\{G\in{[n]\choose \ell}\mid W\subseteq G,\ T\nsubseteq G\right\}.
$$

For each $W\in\mathcal{W}^\prime$, since
$$
\left\{G\in{[n]\choose \ell}\mid W\subseteq G,\ T\nsubseteq G\right\}=\left\{G\in{[n]\choose \ell}\mid W\subseteq G\right\}\setminus\left\{G\in{[n]\choose \ell}\mid W\subseteq G,\ T\subseteq G\right\},
$$
we have
$$
\left|\left\{G\in{[n]\choose \ell}\mid W\subseteq G,\ T\nsubseteq G\right\}\right|={n-t-1\choose\ell-t-1}-{n-t-2\choose \ell-t-2}={n-t-2\choose \ell-t-1}.
$$
 Hence, we have that $|\mathcal{G}^\prime|\leq t(k-t){n-t-2\choose \ell-t-1}$
as desired.     $\qed$

\subsection{$\{\tau_t(\mathcal{F}), \tau_t(\mathcal{G})\}\neq\{t,t+1\}$}

\begin{lemma}\label{s-prer}
Let $n,\ k,\ \ell,\ t$ and $s$ be positive integers with $n\geq k+\ell$, and $G\in{[n]\choose \ell}$ and $S\in{[n]\choose s}$ be two subsets of $[n]$ with $|G\cap S|<t$. Let $\mathcal{F}\subseteq{[n]\choose k}$ be a family satisfying $|G\cap F|\geq t$ for all $F\in\mathcal{F}$. Assume that $|G\cap S|=w$. Then for each $i\in\{1,2,\ldots,t-w\}$, there exists an $(s+i)$-subset $U_i$ with $S\subseteq U_i$  such that  $|\mathcal{F}_S|\leq {\ell-w\choose i} |\mathcal{F}_{U_i}|$. Furthermore,  we have $|\mathcal{F}_{S}|\leq {\ell-w\choose t-w}{n-s-t+w\choose k-s-t+w}$.
\end{lemma}
\proof If $\mathcal{F}_S=\emptyset$, the lemma is clear. Now suppose that $\mathcal{F}_S\neq \emptyset$. For each $i\in\{1,2,\ldots, t-w\}$, write
$$
\mathcal{H}_i=\left\{H\in{G\cup S\choose s+i}\mid S\subseteq H\right\}.
$$
 For each $F\in\mathcal{F}_S$, since $|F\cap G|\geq t$, we have $|F\cap (G\cup S)|=|F\cap G|+|F\cap (S\setminus G)|\geq t+s-w.$ It follows that there exists $H\in\mathcal{H}_i$ such that $H\subseteq F.$ Therefore, $\mathcal{F}_S=\cup_{H\in\mathcal{H}_i}\mathcal{F}_{H}$. Let $U_i$ be an element in $\mathcal{H}_{i}$ such that $|\mathcal{F}_{U_i}|=\max\{|\mathcal{F}_{H}|\mid  H\in\mathcal{H}_i\}$. Observe that $|G\cup S|=\ell+s-w$ and $|\mathcal{H}_i|={\ell-w\choose i}$, and then the former part of this lemma holds. Setting $i=t-w$, then $|U_{t-w}|=s+t-w$, and the latter part of this lemma holds due to $|\mathcal{F}_{U_{t-w}}|\leq {n-s-t+w\choose k-s-t+w}$.   $\qed$
\begin{lemma} \label{s-upper-L}
	Let $n,\ k,\ \ell,\ t,\ w, \ m$ and $s$ be non-negative integers with $1\leq t\leq\ell$.
	\begin{enumerate}[{\rm(i)}]
		\item If $n\geq2(k-t+1)(\ell-t+1)+t+1$ and $t\leq s<k$, then the function ${\ell-w\choose t-w}{n-s-t+w\choose k-s-t+w}$ is increasing as $w\in\{\max\{0,s+t-k\},\ldots,t-1\}$ increases.
		\item If $n\geq (t+1)^2(k-t+1)(\ell-t+1)+t+1$, then the function ${k}^{m-t-2}(k-t+1)^2{m\choose t}{n-m\choose \ell-m}$ is decreasing as $m\in\{t, t+1,\ldots,\ell\}$ increases, and ${k}^{m-t-2}(k-t+1)^2{m\choose t}{n-m\choose \ell-m}<(t+1)(k-t+1){n-t-1\choose \ell-t-1}$ if $t+2\leq m\leq \ell.$
	\end{enumerate}	
\end{lemma}
\proof (i)\quad Let
$$
g(w)={\ell-w\choose t-w}{n-s-t+w\choose k-s-t+w}
$$
for $w\in\{\max\{0,s+t-k\},\ldots,t-1\}$. Since $n\geq 2(k-t+1)(\ell-t+1)+t+1$ and $s\geq t\geq w+1$, we have
\begin{align*}
&(t-w)(n-s-t+w+1)-(\ell-w)(k-s-t+w+1)\\
=& (t-w)(n-k)-(\ell-t)(k-s-t+w+1)\\
\geq& (n-k)-(\ell-t)(k-t)>0
\end{align*}
and
\begin{align*}
\frac{g(w+1)}{g(w)}=&\frac{(t-w)(n-s-t+w+1)}{(\ell-w)(k-s-t+w+1)}\geq 1
\end{align*}
for each $w\in\{\max\{0,s+t-k\},\ldots,t-2\}$. That is, the function $g(w)$ is increasing as $w\in\{\max\{0,s+t-k\},\ldots,t-1\}$ increases.

(ii)\quad For each $m\in\{t,t+1,\ldots,\ell\}$, set
$$
f_2(m,k,\ell,n,t)={k}^{m-t-2}(k-t+1)^2{m\choose t}{n-m\choose \ell-m}.
$$
For each $m\in\{t,t+1,\ldots,\ell-1\}$, since $n\geq(t+1)^2(k-t+1)(\ell-t+1)+t+1$ and $(t+1)(k-t+1)=t(k-t)+k+1>k$, we have
\begin{align*}
\frac{f_2(m+1,k,\ell,n,t)}{f_2(m,k,\ell,n,t)}=&\frac{k(m+1)(\ell-m)}{(m-t+1)(n-m)}=k\left(1+\frac{t}{m-t+1}\right)\left(1-\frac{n-\ell}{n-m}\right)\\
\leq& k(1+t)\left(1-\frac{n-\ell}{n-t}\right)=\frac{k(t+1)(\ell-t)}{n-t}<1,
\end{align*}
implying that the former part of (ii) holds. For $t+2\leq m\leq \ell,$ observe that
$$
\frac{(t+1)(k-t+1){n-t-1\choose \ell-t-1}}{f_2(m,k,\ell,n,t)}>\frac{(t+1)(k-t+1){n-t-1\choose \ell-t-1}}{f_2(t+1,k,\ell,n,t)} =\frac{k}{k-t+1}\geq1,
$$
implying that the latter part of (ii) holds.   $\qed$
\begin{lemma}\label{s-upper-F}
Let $n,$ $k,$ $\ell$ and $t$ be positive integers with $n\geq 2(k-t+1)(\ell-t+1)+t+1$, and $\mathcal{F}\subseteq {[n]\choose k}$ and $\mathcal{G}\subseteq {[n]\choose \ell}$ be maximal cross $t$-intersecting families. Assume that $\tau_{t}(\mathcal{F})=m_f$  and $\tau_{t}(\mathcal{G})=m_g$. Then
	$$
	|\mathcal{F}|\leq{m_f\choose t}{n-t\choose k-t}.
	$$
Furthermore, the following hold.
\begin{enumerate}[{\rm(i)}]
	\item If $m_g=t+1$, then
	$$
	|\mathcal{F}|\leq(\ell-t+1){m_f\choose t}{n-t-1\choose k-t-1}.
	$$
	\item If $m_g\geq t+2$, then
	$$
	|\mathcal{F}|\leq {\ell}^{m_g-t-2}(\ell-t+1)^2{m_f\choose t}{n-m_g\choose k-m_g}.
	$$
\end{enumerate}
\end{lemma}
\proof Let $T$ be a $t$-cover of $\mathcal{F}$ with  size $m_f.$ Then for every $F\in\mathcal{F},$ we have $|F\cap T|\geq t$ and there exists $H\in{T\choose t}$ such that $H\subseteq F,$ implying that $\mathcal{F}\subseteq \cup_{H\in{T\choose t}}\mathcal{F}_{H}$. Let  $H_1$ be an element in ${T\choose t}$ such that $|\mathcal{F}_{H_1}|=\max\{|\mathcal{F}_{H}|\mid H\in{T\choose t}\}$. We have
\begin{align}\label{v-upp-1}
	|\mathcal{F}|\leq{m_f\choose t}|\mathcal{F}_{H_1}|.
\end{align}
Since $|\mathcal{F}_{H_1}|\leq {n-t\choose k-t},$ we have $|\mathcal{F}|\leq{m_f\choose t}{n-t\choose k-t}$.

\

(i)\quad Suppose that $m_g=t+1$. It follows from $|H_1|=t$ and $m_g>t$ that there exists $G_1\in\mathcal{G}$ such that $|G_1\cap H_1|=w<t$. By Lemmas~\ref{s-prer} and \ref{s-upper-L} (i), we have $$|\mathcal{F}_{H_1}|\leq{\ell-w\choose t-w}{n-2t+w\choose k-2t+w}\leq(\ell-t+1){n-t-1\choose k-t-1},$$ implying that (i) holds.

\

(ii)\quad Suppose that $m_{g}\geq t+2$. Firstly we claim that there exists $H^\prime\in{[n]\choose m_g-2}$ such that
\begin{align}\label{s-upp-2}
|\mathcal{F}|\leq\ell^{m_g-t-2}{m_f\choose t}|\mathcal{F}_{H^\prime}|.	
\end{align}
If $m_g=t+2$, it is clear that (\ref{s-upp-2}) holds by setting $H^\prime=H_1$. If $m_g\geq t+3$, using the former part of Lemma~\ref{s-prer} repeatedly, then there exist
$$
H_2\in{[n]\choose t+1},\ H_3\in{[n]\choose t+2},\ldots,H_{m_g-t-1}\in{[n]\choose m_g-2}
$$
such that $H_i\subseteq H_{i+1}$ and $|\mathcal{F}_{H_{i}}|\leq \ell|\mathcal{F}_{H_{i+1}}|$ for each $i\in\{1,2,\ldots,m_g-t-2\}$, which implies that (\ref*{s-upp-2}) holds by setting $H^\prime=H_{m_g-t-1}$. Therefore, the claim holds.

Since $\tau_t(\mathcal{G})>m_g-2=|H^\prime|,$ there exists $G_2\in\mathcal{G}$ such that $|H^\prime\cap G_2|<t.$ If $|H^\prime\cap G_2|=w\leq t-2$, by Lemmas~\ref{s-prer} and \ref{s-upper-L} (i), we have
\begin{align}\label{s-upp-3}
|\mathcal{F}_{H^\prime}|\leq {\ell-w\choose t-w}{n-m_g+2-t+w\choose k-m_g+2-t+w}\leq {\ell-t+2\choose 2}{n-m_g\choose k-m_g}.
\end{align}
Suppose that $|H^\prime\cap G_2|=t-1$. By Lemma~\ref{s-prer}, there exists an $(m_g-1)$-subset $H^{\prime\prime}$ such that $|\mathcal{F}_{H^\prime}|\leq(\ell-t+1)|\mathcal{F}_{H^{\prime\prime}}|$. Since $\tau_t(\mathcal{G})>m_g-1=|H^{\prime\prime}|$,  there exists $G_3\in\mathcal{G}$ such that $|H^{\prime\prime}\cap G_3|<t.$ If $|H^{\prime\prime}\cap G_3|=w\leq t-2$, by Lemma~\ref{s-prer} and \ref{s-upper-L} (i) again, we have
$$
|\mathcal{F}_{H^{\prime\prime}}|\leq {\ell-w\choose t-w}{n-m_g+1-t+w\choose k-m_g+1-t+w}\leq {\ell-t+2\choose 2}{n-m_g-1\choose k-m_g-1},
$$
implying that
\begin{align}\label{s-upp-4}
|\mathcal{F}_{H^\prime}|\leq (\ell-t+1)|\mathcal{F}_{H^{\prime\prime}}|\leq (\ell-t+1){\ell-t+2\choose 2}{n-m_g-1\choose k-m_g-1}.
\end{align}
If $|H^\prime\cap G_3|=t-1$, by Lemma~\ref{s-prer}, then $|\mathcal{F}_{H^{\prime\prime}}|\leq(\ell-t+1){n-m_g\choose k-m_g},$ implying that
\begin{align}\label{s-upp-5}
|\mathcal{F}_{H^\prime}|\leq (\ell-t+1)|\mathcal{F}_{H^{\prime\prime}}|\leq (\ell-t+1)^2{n-m_g\choose k-m_g}.
\end{align}

By $n\geq 2(k-t-1)(\ell+1-t)+t+1$ and $m_g>t+1$  we have
$$
\frac{n-m_{g}}{k-m_{g}}>\frac{n-t-1}{k-t-1}\geq (\ell-t+1)\geq\frac{\ell-t+2}{2},
$$
implying that
$$
(\ell-t+1)^2{n-m_g\choose k-m_g}\geq{\ell-t+2\choose 2}{n-m_g\choose k-m_g}\geq (\ell-t+1){\ell-t+2\choose 2}{n-m_g-1\choose k-m_g-1}.
$$
This together with (\ref{s-upp-2}), (\ref{s-upp-3}), (\ref{s-upp-4}) and (\ref{s-upp-5}) yields (ii) holds. $\qed$
\begin{cor}\label{s-v-cor}
Let $n$, $k$, $\ell$ and $t$ be positive integers with $n\geq (t+1)^2(k-t+1)(\ell-t+1)+t+1$, and $\mathcal{F}\subseteq{[n]\choose k}$ and $\mathcal{G}\subseteq {[n]\choose \ell}$ be maximal cross $t$-intersecting families. Assume that $\tau_t(\mathcal{F})=m_f$ and $\tau_t(\mathcal{G})=m_g$ with $m_f\leq m_g$.
\begin{enumerate}[{\rm (i)}]
\item If $m_f=t$ and $m_g\geq t+2$, then $
|\mathcal{F}||\mathcal{G}|\leq g_4(\ell,k,n,t).$
\item If $m_f\geq t+1,$ then $|\mathcal{F}||\mathcal{G}|<g_5(k,\ell,n,t)$.
\end{enumerate}
\end{cor}
\proof Applying Lemma~\ref{s-upper-F} to $\mathcal{F}$ and $\mathcal{G}$ respectively,  we have
\begin{align}\label{s-upper-F1}
	|\mathcal{F}|\leq\left\{
	\begin{array}{ll}
		(\ell-t+1){m_f\choose t}{n-t-1\choose k-t-1},&\mbox{if}\ m_g=t+1,\vspace{0.1cm}\\
		{\ell}^{m_g-t-2}(\ell-t+1)^2{m_f\choose t}{n-m_g\choose k-m_g},&\mbox{if}\ m_g\geq t+2,
	\end{array}
	\right.
\end{align}
and
\begin{align}\label{s-upper-F2}
	|\mathcal{G}|\leq\left\{
	\begin{array}{ll}
		{m_g\choose t}{n-t\choose \ell-t},&\mbox{if}\ m_f=t,\vspace{0.1cm}\\
		(k-t+1){m_g\choose t}{n-t-1\choose \ell-t-1},&\mbox{if}\ m_f=t+1,\vspace{0.1cm}\\
		{k}^{m_f-t-2}(k-t+1)^2{m_g\choose t}{n-m_f\choose \ell-m_f},&\mbox{if}\ m_f\geq t+2.
	\end{array}
	\right.
\end{align}

(i)\quad If $m_f=t$ and $m_g\geq t+2$, by (\ref{s-upper-F1}), (\ref{s-upper-F2}) and Lemma~\ref{s-upper-L} (ii), then
\begin{align*}
|\mathcal{F}||\mathcal{G}|\leq&{n-t\choose \ell-t}\cdot{\ell}^{m_g-t-2}(\ell-t+1)^2{m_g\choose t}{n-m_g\choose k-m_g}\\
\leq& {n-t\choose \ell-t}\cdot(\ell-t+1)^2{t+2\choose t}{n-t-2\choose k-t-2},
\end{align*}
implying that (i) holds.

(ii)\quad By (\ref{s-upper-F1}), (\ref{s-upper-F2}) and Lemma~\ref{s-upper-L} (ii) again, if $m_f=t+1$ and $m_g=t+1$, then
\begin{align*}
|\mathcal{F}||\mathcal{G}|\leq(t+1)^2(k-t+1)(\ell-t+1){n-t-1\choose \ell-t-1}{n-t-1\choose k-t-1};
\end{align*}
if $m_f=t+1$ and $m_g\geq t+2$, then
\begin{align*}
|\mathcal{F}||\mathcal{G}|\leq&(t+1)(k-t+1){n-t-1\choose \ell-t-1}\cdot{\ell}^{m_g-t-2}(\ell-t+1)^2{m_g\choose t}{n-m_g\choose k-m_g}\\
\leq& (t+1)(k-t+1){n-t-1\choose \ell-t-1}\cdot (t+1)(\ell-t+1){n-t-1\choose k-t-1};
\end{align*}
if $m_g\geq m_f\geq t+2$, then
\begin{align*}
|\mathcal{F}||\mathcal{G}|\leq &{k}^{m_f-t-2}(k-t+1)^2{m_f\choose t}{n-m_f\choose \ell-m_f}\cdot{\ell}^{m_g-t-2}(\ell-t+1)^2{m_g\choose t}{n-m_g\choose k-m_g}\\
\leq & (t+1)(k-t+1){n-t-1\choose \ell-t-1}\cdot (t+1)(\ell-t+1){n-t-1\choose k-t-1}.
\end{align*}
Hence $|\mathcal{F}||\mathcal{G}|<g_5(k,\ell,n,t)$ holds.  $\qed$

\subsection{The proof of Theorem~\ref{s-main-1}}
In this part, we will prove Theorem~\ref{s-main-1}. Let $n$, $k_1$, $k_2$ and $t$ be positive integers satisfying $n\geq\max\{(t+1),(k_2-t)\}\cdot(t+1)(k_1-t+1)(k_2-t+1)+t+1$, $k_1\geq k_2\geq t+1$ and $(k_1,k_2,t)\neq (2,2,1),(3,2,1),(4,2,1)$ or $(4,4,2)$. Assume that $\mathcal{F}_1\subseteq{[n]\choose k_1}$ and $\mathcal{F}_2\subseteq{[n]\choose k_2}$ are maximal non-trivial cross $t$-intersecting families.

By Lemma~\ref{s-FGtt+1}, note that $|\mathcal{F}_1||\mathcal{F}_2|=g_1(k_1,k_2,n,t)$ if $\mathcal{F}_1$ and $\mathcal{F}_2$ are a pair of families given in ${\rm (ia)}$ or ${\rm (ib)},$ and  $|\mathcal{F}_1||\mathcal{F}_2|=g_2(k_1,k_2,n,t)$ if $\mathcal{F}_1$ and $\mathcal{F}_2$ are a pair of families given in ${\rm (iia)}$ or ${\rm (iib)}.$ From Lemmas~\ref{s-ineq-lem-2}, \ref{s-ineq-lem-11}, \ref{s-ineq-lem-3} and  \ref{s-ineq-lem-4}, we have
\begin{align}\label{s-v-ineq-10}
g_1(k_1,k_2,n,t)\left\{\begin{array}{ll}
>g_2(k_1,k_2,n,t),& \mbox{if}\ k_2\geq 2t+1,\\
<g_2(k_1,k_2,n,t),& \mbox{if}\  k_2\leq 2t.
\end{array}\right.
\end{align}
Now suppose that $\mathcal{F}_1$ and $\mathcal{F}_2$ are maximal non-trivial cross $t$-intersecting families, which are neither a pair of families given in ${\rm(ia)}$ and ${\rm(ib)}$, nor a pair of families given in ${\rm(iia)}$ and ${\rm(iib)}$. To prove the theorem, we only need to show that
\begin{align*}
|\mathcal{F}_1||\mathcal{F}_2|<\left\{\begin{array}{ll}
g_1(k_1,k_2,n,t),& \mbox{if}\ k_2\geq 2t+1,\\
g_2(k_1,k_2,n,t),& \mbox{if}\  k_2\leq 2t.
\end{array}\right.
\end{align*}
By (\ref{s-v-ineq-10}), it suffices to prove that
\begin{align}\label{s-F1F2upp}
|\mathcal{F}_1||\mathcal{F}_2|<g_1(k_1,k_2,n,t)\quad \mbox{or}\quad  |\mathcal{F}_1||\mathcal{F}_2|<g_2(k_1,k_2,n,t).
\end{align}
 We divide our proof into the following five cases.

\medskip

\noindent\textbf{Case 1.} $\tau_t(\mathcal{F}_1)=\tau_t(\mathcal{F}_2)=t$.

Assume that $X_1$ and $X_2$ are $t$-covers of $\mathcal{F}_1$ and $\mathcal{F}_2$ with size $t$, respectively. By Lemma~\ref{s-tt}, we have $|X_1\cap X_2|\geq t$, implying that $X_1=X_2$. Then $\mathcal{F}_1$ and $\mathcal{F}_2$ are trivial cross $t$-intersecting families, a contradiction.

\medskip

\noindent\textbf{Case 2.} $\tau_t(\mathcal{F}_1)=t$ and $\tau_t(\mathcal{F}_2)=t+1$.

Assume that $X$ is a $t$-cover of $\mathcal{F}_1$ with size $t$, and $T$ is a $t$-cover of $\mathcal{F}_2$ with size $t+1$. By Lemma~\ref{s-FGtt+1}, we consider the following subcases.

\noindent\textbf{Case 2.1.}\quad There exists a $(k_2+1)$-subsset $M$ of $[n]$ such that $F_2\cup X=M$ for each $F_2\in\mathcal{F}_2\setminus(\mathcal{F}_2)_X$.

By Lemma~\ref{s-FGtt+1} (ii), we have
\begin{align*}
\mathcal{F}_1=&\left\{F_1\in{[n]\choose k_1}\mid X\subseteq F_1,\ |F_1\cap M|\geq t+1\right\},\\
\mathcal{F}_2=&\left\{F_2\in{[n]\choose k_2} \mid X\subseteq F_2\right\}\cup{M\choose k_2},
\end{align*}
 and $|\mathcal{F}_1||\mathcal{F}_2|=g_1(k_1,k_2,n,t)$. Since $\mathcal{F}_1$ and $\mathcal{F}_2$ are not a pair of families given in ${\rm (ia)}$ or ${\rm (ib)},$ we have $t+1\leq k_2\leq 2t$. It follows from Lemmas~\ref{s-ineq-lem-11}, \ref{s-ineq-lem-3} and \ref{s-ineq-lem-4} that
$|\mathcal{F}_1||\mathcal{F}_2|<g_2(k_1,k_2,n,t).$ Therefore, (\ref{s-F1F2upp}) holds.

\noindent\textbf{Case 2.2.}\quad $T\subseteq F_1$ for each $F_1\in\mathcal{F}_1.$

By Lemma~\ref{s-FGtt+1} (iii), we have $\mathcal{F}_1=\left\{F_1\in{[n]\choose k_1}\mid T\subseteq F_1\right\}$, $\mathcal{F}_2=\left\{F_2\in{[n]\choose k_2}\mid |F_2\cap T|\geq t\right\}$ and $|\mathcal{F}_1||\mathcal{F}_2|=g_2(k_1,k_2,n,t)$. Since $\mathcal{F}_1$ and $\mathcal{F}_2$ are not a pair of families given in ${\rm (iia)}$ or ${\rm (iib)}$, we have $k_2\geq 2t+1$. It follows from Lemma~\ref{s-ineq-lem-2} that
$|\mathcal{F}_1||\mathcal{F}_2|<g_1(k_1,k_2,n,t),$
and so (\ref{s-F1F2upp}) holds.

\noindent\textbf{Case 2.3.}\quad There exist $F_{2,1},\ F_{2,2}\in\mathcal{F}_2\setminus(\mathcal{F}_2)_X$ such that $F_{2,1}\cup X\neq F_{2,2}\cup X$, and there exists $F_1\in\mathcal{F}_1$ such that $T\nsubseteq F_1$.

It follows from Lemmas~\ref{s-FGtt+1} (iv) and \ref{s-ineq-lem-5}  that
$$
|\mathcal{F}_1||\mathcal{F}_2|\leq g_3(k_1,k_2,n,t)<g_1(k_1,k_2,n,t),
$$
and so (\ref{s-F1F2upp}) holds.

\medskip

\noindent\textbf{Case 3.} $\tau_t(\mathcal{F}_1)=t+1$ and $\tau_t(\mathcal{F}_2)=t$.

Assume that $T$ is a $t$-cover of $\mathcal{F}_1$ with size $t+1$, and $X$ is a $t$-cover of $\mathcal{F}_2$ with size $t$. By Lemma~\ref{s-FGtt+1}, we consider the following subcases.

\noindent\textbf{Case 3.1.}\quad There exists a $(k_1+1)$-subset $M$ of $[n]$ such that $F_1\cup X=M$ for each $F_1\in\mathcal{F}_1\setminus(\mathcal{F}_1)_X$.

By Lemma~\ref{s-FGtt+1} (ii), we have
\begin{align*}
\mathcal{F}_1=&\left\{F_1\in{[n]\choose k_1}\mid X\subseteq F_1\right\}\cup{M\choose k_1},\\
\mathcal{F}_2=&\left\{F_2\in{[n]\choose k_2}\mid X\subseteq F_2, \ |F_2\cap M|\geq t+1\right\},
\end{align*}
and $|\mathcal{F}_1||\mathcal{F}_2|=g_1(k_2,k_1,n,t)$. If $k_1>k_2$,  then $|\mathcal{F}_1||\mathcal{F}_2|<g_1(k_1,k_2,n,t)$ from Lemma~\ref{s-ineq-lem-6}. Now assume that $k_1=k_2$. Since $\mathcal{F}_1$ and $\mathcal{F}_2$ are not a pair of families given in ${\rm (ib)},$ we have $k_2\leq 2t$, implying that $|\mathcal{F}_1||\mathcal{F}_2|=g_1(k_2,k_1,n,t)=g_1(k_1,k_2,n,t)<g_2(k_1,k_2,n,t)$ from (\ref{s-v-ineq-10}). Therefore (\ref{s-F1F2upp}) holds.

\noindent\textbf{Case 3.2.}\quad $T\subseteq F_2$ for each $F_2\in\mathcal{F}_2.$

By Lemma~\ref{s-FGtt+1} (iii), we have $\mathcal{F}_1=\left\{F_1\in{[n]\choose k_1}\mid |F_1\cap T|\geq t\right\}$, $\mathcal{F}_2=\left\{F_2\in{[n]\choose k_2}\mid T\subseteq F_2\right\}$ and $|\mathcal{F}_1||\mathcal{F}_2|=g_2(k_2,k_1,n,t)$. If $k_1>k_2$, then $|\mathcal{F}_1||\mathcal{F}_2|<g_2(k_1,k_2,n,t)$ from Lemma~\ref{s-ineq-lem-7}. Now assume that $k_1=k_2$. Since $\mathcal{F}_1$ and $\mathcal{F}_2$ are not a pair of families given in ${\rm (iib)},$ we have $k_2\geq 2t+1$, implying that $|\mathcal{F}_1||\mathcal{F}_2|=g_2(k_2,k_1,n,t)=g_2(k_1,k_2,n,t)<g_1(k_1,k_2,n,t)$ by (\ref{s-v-ineq-10}). Therefore (\ref{s-F1F2upp}) holds.

\noindent\textbf{Case 3.3.}\quad There exist $F_{1,1},\ F_{1,2}\in\mathcal{F}_1\setminus(\mathcal{F}_1)_X$ such that $F_{1,1}\cup X\neq F_{1,2}\cup X$, and there exists $F_2\in\mathcal{F}_2$ such that $T\nsubseteq F_2$.

It follows from Lemmas~\ref{s-FGtt+1} (iv) and \ref{s-ineq-lem-5} that
$$
|\mathcal{F}_1||\mathcal{F}_2|\leq g_3(k_2,k_1,n,t)<g_1(k_1,k_2,n,t),
$$
and hence (\ref{s-F1F2upp}) holds.

\medskip

\noindent\textbf{Case 4.} $\tau_t(\mathcal{F}_1)=t$ and $\tau_t(\mathcal{F}_2)\geq t+2$, or $\tau_t(\mathcal{F}_1)\geq t+2$ and $\tau_t(\mathcal{F}_2)=t$.

Assume that $\tau_t(\mathcal{F}_1)=t$ and $\tau_t(\mathcal{F}_2)\geq t+2$. By Corollary~\ref{s-v-cor} (i) and Lemma~\ref{s-ineq-lem-8}, we have that
$$
|\mathcal{F}_1||\mathcal{F}_2|\leq g_4(k_2,k_1,n,t)<g_1(k_1,k_2,n,t),
$$
and (\ref{s-F1F2upp}) holds.

Assume that  $\tau_t(\mathcal{F}_1)\geq t+2$ and $\tau_t(\mathcal{F}_2)=t$. By Corollary~\ref{s-v-cor}  (i) and Lemma~\ref{s-ineq-lem-8} again, we have that
$$
|\mathcal{F}_1||\mathcal{F}_2|\leq g_4(k_1,k_2,n,t)<g_1(k_1,k_2,n,t),
$$
and (\ref{s-F1F2upp}) holds.

\medskip

\noindent\textbf{Case 5.} $\tau_t(\mathcal{F}_1)\geq t+1$ and $\tau_t(\mathcal{F}_2)\geq t+1$.			

From Corollary~\ref{s-v-cor} (ii) and Lemma~\ref{s-ineq-lem-10}, we obtain
$$
|\mathcal{F}_1||\mathcal{F}_2|<g_5(k_1,k_2,n,t)<g_1(k_1,k_2,n,t),
$$
and (\ref{s-F1F2upp}) holds.    $\qed$

\section{Some inequalities}

In this section, we prove some inequalities used in the proof of Theorem~\ref{s-main-1}. In the following, we always assume that $n$, $k_1$, $k_2$ and $t$ are positive integers with $k_1\geq k_2\geq t+1$ and $n\geq \max\{t+1,\ k_2-t\}\cdot(t+1)(k_1-t+1)(k_2-t+1)+t+1$. The following Pascal's formula is used frequently:
{\small\begin{align*}
{a\choose b}={a-1\choose b}+{a-1\choose b-1},
\end{align*}}
where $a$ and $b$ are positive integers with $a\geq b.$

Let $g_i(k,\ell,n,t)\ (1\leq i\leq 5)$ be as in (\ref{s-eq-1})--(\ref{s-eq-5}), and write
{\small\begin{align}
g_6(k,\ell,n,t)=&\left((\ell-t+1){n-t-1\choose k-t-1}-{\ell-t+1\choose 2}{n-t-2\choose k-t-2}\right){n-t\choose \ell-t}.\label{s-eq-6}
\end{align}}%
For $i\in\{1,2,\ldots,6\}$, set
{\small\begin{align*}
\tilde{g}_i(k,\ell,n,t)=g_i(k,\ell,n,t){n-t-1\choose k-t-1}^{-1}{n-t-1\choose \ell-t-1}^{-1}.
\end{align*}}%
Observe that
{\small\begin{align}\label{s-tg2}
\tilde{g}_2(k_1,k_2,n,t)=\frac{(t+1)(n-k_2)}{k_2-t}+1\quad\mbox{and}\quad\tilde{g}_2(k_2,k_1,n,t)=\frac{(t+1)(n-k_1)}{k_1-t}+1,
\end{align}}%
{\small\begin{align}\label{s-tg6}
\tilde{g}_6(k_1,k_2,n,t)=&(k_2-t+1)\left(\frac{n-t}{k_2-t}-\frac{k_1-t}{2}+\frac{n-k_1}{2(n-t-1)}\right)\nonumber\\
>&(k_2-t+1)\left(\frac{n-t}{k_2-t}-\frac{k_1-t}{2}\right).
\end{align}}%
\begin{lemma}\label{s-ineq-lem-1}
We have $g_6(k_1,k_2,n,t)<g_1(k_1,k_2,n,t)<(k_2-t+1){n-t-1\choose k_1-t-1}\left({n-t\choose k_2-t}+t\right).$
\end{lemma}
\proof By (\ref{s-eq-1}) and (\ref{s-eq-6}), it suffices to prove that
{\small\begin{align}
&{n-t\choose k_1-t}-{n-k_2-1\choose k_1-t}<(k_2-t+1){n-t-1\choose k_1-t-1},\label{s-ineq-2}\\
&(k_2-t+1){n-t-1\choose k_1-t-1}-{k_2-t+1\choose 2}{n-t-2\choose k_1-t-2}<{n-t\choose k_1-t}-{n-k_2-1\choose k_1-t}\label{s-ineq-1}.
\end{align}}%

Using Pascal's formula repeatedly, we obtain
 {\small\begin{align*}
 {n-t\choose k_1-t}-{n-k_2-1\choose k_1-t}=\sum_{i=0}^{k_2-t}{n-t-1-i\choose k_1-t-1}\leq (k_2-t+1){n-t-1\choose k_1-t-1}
 \end{align*} }%
and
 {\small\begin{align*}
(k_2-t+1){n-t-1\choose k_1-t-1}=&{n-t-1\choose k_1-t-1}+\sum_{i=1}^{k_2-t}\left(\sum_{j=1}^{i}{n-t-1-j\choose k_1-t-2}+{n-t-1-i\choose k_1-t-1}\right)\\
=&\sum_{i=0}^{k_2-t}{n-t-1-i\choose k_1-t-1}+\sum_{i=1}^{k_2-t}\sum_{j=1}^{i}{n-t-1-j\choose k_1-t-2}\\
\leq&{n-t\choose k_1-t}-{n-k_2-1\choose k_1-t}+{k_2-t+1\choose 2}{n-t-2\choose k_1-t-2},
 \end{align*} }%
implying that (\ref{s-ineq-2}) and (\ref{s-ineq-1}) hold.    $\qed$
\begin{lemma}\label{s-ineq-lem-2}
If $k_2\geq 2t+1,$ then $g_1(k_1,k_2,n,t)>g_2(k_1,k_2,n,t)$.
\end{lemma}
\proof By (\ref{s-tg2}) and (\ref{s-tg6}), we have
{\small\begin{align*}
&\tilde{g}_6(k_1,k_2,n,t)-\tilde{g}_2(k_1,k_2,n,t)\\
>&(k_2-t+1)\left(\frac{n-t}{k_2-t}-\frac{k_1-t}{2}\right)-\frac{(t+1)(n-k_2)}{k_2-t}-1\\
=&\frac{(k_2-t+1)(n-t)}{k_2-t}-\frac{(k_2-t+1)(k_1-t)}{2}-\frac{(t+1)(n-t)}{k_2-t}+t\\
>&\frac{(k_2-2t)(n-t)}{k_2-t}-\frac{(k_2-t+1)(k_1-t)}{2}>0
\end{align*}}%
due to $k_2\geq 2t+1$ and $n\geq (k_2-t)(t+1)(k_1-t+1)(k_2-t+1)+t+1$. Then $g_6(k_1,k_2,n,t)>g_2(k_1,k_2,n,t)$, implying that the required lemma holds by Lemma~\ref{s-ineq-lem-1}. $\qed$
\begin{lemma}\label{s-ineq-lem-11}
If $k_2=2t$, $t\geq 2$ and $(k_1,k_2)\neq(4,4)$, then $g_2(k_1,k_2,n,t)>g_1(k_1,k_2,n,t)$.
\end{lemma}%
\proof First, observe that
{\small\begin{align*}
\frac{g_2(k_1,2t,n,t)\cdot g_1(k_1-1,2t,n,t)}{g_1(k_1,2t,n,t)\cdot g_2(k_1-1,2t,n,t)}=&\frac{{n-t-1\choose k_1-t-1}\left({n-t\choose k_1-t-1}-{n-2t-1\choose k_1-t-1}\right)}{{n-t-1\choose k_1-t-2}\left({n-t\choose k_1-t}-{n-2t-1\choose k_1-t}\right)}\\
=&\frac{(n-k_1+1)\left({n-t\choose k_1-t-1}-{n-2t-1\choose k_1-t-1}\right)}{(k_1-t-1)\left({n-t\choose k_1-t}-{n-2t-1\choose k_1-t}\right)}.
\end{align*}}%
Since $(n-k_1+1){n-t\choose k_1-t-1}=(k_1-t){n-t\choose k_1-t}$ and $(n-k_1+1){n-2t-1\choose k_1-t-1}=(n-k_1-t){n-2t-1\choose k_1-t-1}+(t+1){n-2t-1\choose k_1-t-1}=(k_1-t){n-2t-1\choose k_1-t}+(t+1){n-2t-1\choose k_1-t-1}$, we have
{\small\begin{align*}
&(n-k_1+1)\left({n-t\choose k_1-t-1}-{n-2t-1\choose k_1-t-1}\right)-(k_1-t-1)\left({n-t\choose k_1-t}-{n-2t-1\choose k_1-t}\right)\\
=&{n-t\choose k_1-t}-{n-2t-1\choose k_1-t}-(t+1){n-2t-1\choose k_1-t-1}\\
=&\sum_{i=0}^t{n-t-1-i\choose k_1-t-1}-(t+1){n-2t-1\choose k_1-t-1}>0
\end{align*}}%
from  ${n-t\choose k_1-t}=\sum_{i=0}^t{n-t-1-i\choose k_1-t-1}+{n-2t-1\choose k_1-t}.$ It follows that $\frac{g_2(k_1,2t,n,t)}{g_1(k_1,2t,n,t)}> \frac{g_2(k_1-1,2t,n,t)}{g_1(k_1-1,2t,n,t)}$ for $2t+1\leq k_1\leq \frac{n-t-1}{(t^2+2t+1)(t+1)}+t-1$. In order to prove the lemma, by induction, it suffices to prove that $g_2(5,4,n,2)>g_1(5,4,n,2)$ if $t=2$, and $g_2(2t,2t,n,t)>g_1(2t,2t,n,t)$ if $t\geq 3.$

\textbf{Step 1.} Show that $g_2(5,4,n,2)>g_1(5,4,n,2)$.

Since ${n-2\choose 3}-{n-5\choose 3}={n-3\choose 2}+{n-4\choose 2}+{n-5\choose 2}<3{n-3\choose 2}$ and $3{n-3\choose 2}-\left({n-2\choose 3}-{n-5\choose 3}\right)=2{n-3\choose 2}-{n-4\choose 2}-{n-5\choose 2}=3n-13$, by (\ref{s-eq-1}) and (\ref{s-eq-2}), we have
{\small\begin{align*}
&(g_2(5,4,n,2)-g_1(5,4,n,2)){n-3\choose 2}^{-1}\\
=&3{n-3\choose 2}+n-3-\left({n-2\choose 3}-{n-5\choose 3}\right)\left({n-2\choose 2}+2\right){n-3\choose 2}^{-1}\\
=&3{n-3\choose 2}+n-3-\left({n-2\choose 3}-{n-5\choose 3}\right)\left({n-3\choose 2}+n-1\right){n-3\choose 2}^{-1}\\
>&3{n-3\choose 2}+n-3-\left({n-2\choose 3}-{n-5\choose 3}\right)-3(n-1)\\
=&n-13>0
\end{align*}}%
from $n\geq (t+1)^2(k_1-t+1)(k_2-t+1)+t+1=111$.

\

\textbf{Step 2.} Show that $g_2(2t,2t,n,t)>g_1(2t,2t,n,t)$ if $t\geq 3.$

  Using Pascal's formula repeatedly, we obtain
 {\small \begin{align*}
  {n-t\choose t}-{n-2t-1\choose t}=\sum_{i=0}^t{n-t-1-i\choose t-1}<(t+1){n-t-1\choose t-1},
   \end{align*}}%
   and
  {\small \begin{align*}
   t{n-t-1\choose t-1}=&\sum_{i=1}^t\left(\sum_{j=1}^i{n-t-1-j\choose t-2}+{n-t-1-i\choose t-1}\right)\\
   =&\sum_{i=1}^t{n-t-1-i\choose t-1}+\sum_{j=1}^t(t+1-j){n-t-1-j\choose t-2}.
   \end{align*}}%
 It follows from (\ref{s-eq-1}) and (\ref{s-eq-2}) that
{\small\begin{align*}
&g_2(2t,2t,n,t)-g_1(2t,2t,n,t)\\
=&{n-t-1\choose t-1}\left(t{n-t-1\choose t}+{n-t\choose t}\right)-\left({n-t\choose t}-{n-2t-1\choose t}\right)\left({n-t\choose t}+t\right)\\
>&{n-t-1\choose t-1}\left(t{n-t-1\choose t}+{n-t\choose t}\right)-\sum_{i=0}^t{n-t-1-i\choose t-1}{n-t\choose t}-t\sum_{i=0}^t{n-t-1-i\choose t-1}\\
=&t{n-t-1\choose t-1}{n-t-1\choose t}-\sum_{i=1}^t{n-t-1-i\choose t-1}\left({n-t-1\choose t}+{n-t-1\choose t-1}\right)-t(t+1){n-t-1\choose t-1}\\
=&\sum_{i=1}^t(t+1-i){n-t-1-i\choose t-2}{n-t-1\choose t}-\sum_{i=1}^t{n-t-1-i\choose t-1}{n-t-1\choose t-1}-t(t+1){n-t-1\choose t-1}.
\end{align*}}%

For each $i\in\{2,\ldots,t-1\}$, observe that
{\begin{align*}
\frac{(t+1-i){n-t-1-i\choose t-2}{n-t-1\choose t}}{{n-t-1-i\choose t-1}{n-t-1\choose t-1}}=\frac{(t+1-i)(t-1)(n-2t)}{t(n-2t-i+1)}>\frac{2(t-1)}{t}>1
\end{align*}}%
which implies that
{\small\begin{align*}
&\sum_{i=1}^t(t+1-i){n-t-1-i\choose t-2}{n-t-1\choose t}-\sum_{i=1}^t{n-t-1-i\choose t-1}{n-t-1\choose t-1}\\
>&\left(t{n-t-2\choose t-2}+{n-2t-1\choose t-2}\right){n-t-1\choose t}-\left({n-t-2\choose t-1}+{n-2t-1\choose t-1}\right){n-t-1\choose t-1}\\
=&{n-t-1\choose t-1}\left((n-2t){n-t-2\choose t-2}+\frac{n-2t}{t}{n-2t-1\choose t-2}-{n-t-2\choose t-1}-{n-2t-1\choose t-1}\right)\\
=&{n-t-1\choose t-1}\left((t-2){n-t-2\choose t-1}-{n-2t-1\choose t-1}+\frac{n-2t}{t}{n-2t-1\choose t-2}\right)\\
>&\frac{n-2t}{t}{n-2t-1\choose t-2}{n-t-1\choose t-1}
\end{align*}}%
due to $(t-2){n-t-2\choose t-1}-{n-2t-1\choose t-1}>0$.

Since $t\geq 3$ and $n\geq(t+1)^2(k_1-t+1)(k_2-t+1)+t+1=(t+1)^4+t+1$, we have $n-2t>n-2t-1\geq(t+1)^4-t>t(t+1)^3$, implying that
{\small\begin{align*}
\frac{n-2t}{t}{n-2t-1\choose t-2}\geq\frac{(n-2t)(n-2t-1)}{t(t-2)}>\frac{t^2(t+1)^6}{t(t-2)}>t(t+1).
\end{align*}}%
Therefore, we have
{\small\begin{align*}
&g_2(2t,2t,n,t)-g_1(2t,2t,n,t)>\frac{n-2t}{t}{n-2t-1\choose t-2}{n-t-1\choose t-1}-t(t+1){n-t-1\choose t-1}>0,
\end{align*}}%
and the required result follows.   $\qed$
\begin{lemma}\label{s-ineq-lem-3}
If $t+2\leq k_2\leq 2t-1,$ then $g_2(k_1,k_2,n,t)>g_1(k_1,k_2,n,t).$
\end{lemma}
\proof Since $n\geq (t+1)^2(k_1-t+1)(k_2-t+1)+t+1$ and $k_1\geq k_2\geq t+2$, we have ${n-t-1\choose k_2-t-1}\geq\frac{n-t-1}{k_2-t-1}\geq(t+1)^2(k_1-t+1)>t(k_2-t+1),$ implying that ${n-t-1\choose k_1-t-1}{n-t-1\choose k_2-t-1}>t(k_2-t+1){n-t-1\choose k_1-t-1}.$ Then, by (\ref{s-eq-2}) and Lemma~\ref{s-ineq-lem-1}, we obtain
{\small\begin{align*}
&g_2(k_1,k_2,n,t)-g_1(k_1,k_2,n,t)\\
>&{n-t-1\choose k_1-t-1}\left((t+1){n-t-1\choose k_2-t}+{n-t-1\choose k_2-t-1}\right)-(k_2-t+1){n-t-1\choose k_1-t-1}\left({n-t\choose k_2-t}+t\right)\\
>&(t+1){n-t-1\choose k_1-t-1}{n-t-1\choose k_2-t}-(k_2-t+1){n-t-1\choose k_1-t-1}{n-t\choose k_2-t}.
\end{align*}}%
It follows from $t+2\leq k_2\leq 2t-1$ and $n-t^2>(t+1)^2(k_1-t+1)(k_2-t+1)-t^2>k_2-t$ that
{\small\begin{align*}
&\left(g_2(k_1,k_2,n,t)-g_1(k_1,k_2,n,t)\right){n-t-1\choose k_1-t-1}^{-1}{n-t-1\choose k_2-t-1}^{-1}\\
=&\frac{(t+1)(n-k_2)}{k_2-t}-\frac{(k_2-t+1)(n-t)}{k_2-t}=\frac{(2t-k_2)n-t^2}{k_2-t}-1>0.
\end{align*}}%
Hence, the required result follows.   $\qed$
\begin{lemma}\label{s-ineq-lem-4}
If $k_2=t+1,$ and $(k_1, t)\neq(2,1),\ (3,1)$ or $(4,1)$, then $g_2(k_1,k_2,n,t)>g_1(k_1,k_2,n,t).$
\end{lemma}
\proof By (\ref{s-eq-1}) and (\ref{s-eq-2}), then
{\small\begin{align*}
g_1(k_1,k_2,n,t)=&n\left({n-t-1\choose k_1-t-1}+{n-t-2\choose k_1-t-1}\right), \\ g_2(k_1,k_2,n,t)=&(t(n-t-2)+n){n-t-1\choose k_1-t-1}.
\end{align*}}%
It follows from $n\geq (t+1)^2(k_1-t+1)(k_2-t+1)+t+1>4(t+1)^2>t(2t+3)$ and $(k_1, t)\neq(2,1),\ (3,1)$ or $(4,1)$ that
{\small\begin{align*}
&(g_2(k_1,k_2,n,t)-g_1(k_1,k_2,n,t))(n-k_1){n-t-2\choose k_1-t-1}^{-1}\\
=&t(n-t-2)(n-t-1)-n(n-k_1)=n((t-1)n+k_1-t(2t+3))+t(t+2)(t+1)>0.
\end{align*}}%
Hence, the lemma holds. $\qed$
\begin{lemma}\label{s-ineq-lem-6}
If $k_1>k_2$, then $g_1(k_1,k_2,n,t)>g_1(k_2,k_1,n,t).$
\end{lemma}
\proof We prove this lemma from the following two cases.

\noindent\textbf{Case 1.} $k_2=t+1$.

By (\ref{s-eq-1}), we have
{\small\begin{align*}
&g_1(k_1,k_2,n,t)-g_1(k_2,k_1,n,t)\\
=&n\left({n-t\choose k_1-t}-{n-t-2\choose k_1-t}\right)-(k_1+1-t)\left({n-t\choose k_1-t}+t\right)\\
=&(n-k_1+t-1){n-t\choose k_1-t}-(n+k_1-t-1){n-t-2\choose k_1-t}+(k_1-t-1){n-t-2\choose k_1-t}-t(k_1+1-t).
\end{align*}}%

Observe that
{\small\begin{align*}
&(n-k_1+t-1)(n-t)(n-t-1)-(n+k_1-t-1)(n-k_1)(n-k_1-1)\\
=&nt(k_1-t+1)-t(t+1)(k_1-t+1)+nk_1(k_1-t-1)-k_1(k_1+1)(k_1-t-1)>0
\end{align*}}%
from $n\geq (t+1)^2(k_1-t+1)(k_2-t+1)+t+1>\max\{t+1, k_1+1\}$, which implies that
{\small\begin{align*}
\frac{(n-k_1+t-1){n-t\choose k_1-t}}{(n+k_1-t-1){n-t-2\choose k_1-t}}=\frac{(n-k_1+t-1)(n-t)(n-t-1)}{(n+k_1-t-1)(n-k_1)(n-k_1-1)}>1.
\end{align*}}%

Since $n\geq 2(t+1)^2(k_1-t+1)+t+1$, we have
{\small\begin{align*}
n-t-2> 2(t^2+2t)(k_1-t+1)\quad\mbox{and}\quad n-t-3\geq 2(t+1)^2(k_1-t),
\end{align*}}%
implying that
{\small\begin{align*}
(k_1-t-1){n-t-2\choose k_1-t}-t(k_1+1-t)\geq\frac{(n-t-2)(n-t-3)}{k_1-t}-t(k_1+1-t)>0
\end{align*}}%
due to $k_1\geq t+2.$ Hence, we obtain $g_1(k_1,k_2,n,t)-g_1(k_2,k_1,n,t)>0$ as desired.

\

\noindent\textbf{Case 2.} $k_2\geq t+2.$

Since
{\small\begin{align*}
\frac{{n-t\choose k_1-t}{n-k_1-1\choose k_2-t}}{{n-k_2-1\choose k_1-t}{n-t\choose k_2-t}}=\frac{(n-t)\cdots(n-k_1+1)\cdot(n-k_1-1)\cdots(n-k_1-k_2+t)}{(n-k_2-1)\cdots(n-k_1-k_2+t)\cdot(n-t)\cdots(n-k_2+1)}=\frac{n-k_2}{n-k_1},
\end{align*}}%
we have ${n-t\choose k_1-t}{n-k_1-1\choose k_2-t}=\frac{n-k_2}{n-k_1}{n-k_2-1\choose k_1-t}{n-t\choose k_2-t}$. By (\ref{s-eq-1}) and $n-t>2(k_1-t)$, we have
{\small\begin{align*}
&g_1(k_1,k_2,n,t)-g_1(k_2,k_1,n,t)\\
=&{n-t\choose k_1-t}{n-k_1-1\choose k_2-t}-{n-k_2-1\choose k_1-t}{n-t\choose k_2-t}+t{n-t\choose k_1-t}-t{n-t\choose k_2-t}\\
&+t{n-k_1-1\choose k_2-t}-t{n-k_2-1\choose k_1-t}\\
>&{n-t\choose k_1-t}{n-k_1-1\choose k_2-t}-{n-k_2-1\choose k_1-t}{n-t\choose k_2-t}-t{n-k_2-1\choose k_1-t}\\
=&\left(\frac{k_1-k_2}{n-k_1}{n-t\choose k_2-t}-t\right){n-k_2-1\choose k_1-t}>0
\end{align*}}%
due to
{\small\begin{align*}
\frac{k_1-k_2}{n-k_1}{n-t\choose k_2-t}\geq\frac{(k_1-k_2)(n-t)(n-t-1)}{(n-k_1)(k_2-t)(k_2-t-1)}> \frac{n-t-1}{(k_2-t)(k_2-t-1)}>t.
\end{align*}}%
Therefore, the required result follows.   $\qed$
\begin{lemma}\label{s-ineq-lem-7}
If $k_1>k_2$, then $g_2(k_1,k_2,n,t)>g_2(k_2,k_1,n,t).$
\end{lemma}
\proof By (\ref{s-tg2}), observe that
{\small\begin{align*}
\tilde{g}_2(k_1,k_2,n,t)-\tilde{g}_2(k_2,k_1,n,t)=&\frac{(t+1)(n-k_2)}{k_2-t}-\frac{(t+1)(n-k_1)}{k_1-t}\\
=&\frac{(t+1)(n-t)(k_1-k_2)}{(k_1-t)(k_2-t)}>0,
\end{align*}}%
implying that the desired result holds. $\qed$
\begin{lemma}\label{s-ineq-lem-5}
We have $g_1(k_1,k_2,n,t)>\max\{g_3(k_1,k_2,n,t),g_3(k_2,k_1,n,t)\}$.
\end{lemma}
\proof It suffices to prove that $g_6(k_1,k_2,n,t)>\max\{g_3(k_1,k_2,n,t),g_3(k_2,k_1,n,t)\}$ by Lemma~\ref{s-ineq-lem-1}.
Set $(k,\ell)=(k_1, k_2)$ or $(k_2,k_1)$. By (\ref{s-eq-3}), observe that
{\small\begin{align*}
\tilde{g}_3(k,\ell, n,t)=&\left(\ell-t+\frac{k-t-1}{n-t-1}\right)\left(\frac{n-t}{\ell-t}+\frac{t(k-t)(n-\ell)}{n-t-1}\right)\\
\leq&\left(\ell-t+\frac{k-t-1}{n-t-1}\right)\left(\frac{n-t}{\ell-t}+t(k-t)\right)\\
=&n-t+t(k-t)(\ell-t)+\frac{(k-t-1)(n-t)}{(n-t-1)(\ell-t)}+\frac{t(k-t)(k-t-1)}{n-t-1}\\
<&n-t+(k-t)(t\ell-t^2+1)+\frac{t(k-t)(k-t-1)}{n-t-1}
\end{align*}}%
due to $\ell\geq t+1$ and $\frac{(k-t-1)(n-t)}{(n-t-1)(\ell-t)}<\frac{k-t}{\ell-t}\leq k-t$.

Since $n\geq(t+1)(k_2-t)(k_1-t+1)(k_2-t+1)+t+1$ and $k_2\geq t+1$, by (\ref{s-tg6}), we have
{\small\begin{align*}
&\tilde{g}_6(k_1,k_2, n,t)-\tilde{g}_3(k_1,k_2, n,t)\\
>&\frac{n-t}{k_2-t}-\frac{(k_1-t)(k_2-t+1)}{2}-(k_1-t)(tk_2-t^2+1)-\frac{t(k_1-t)(k_1-t-1)}{n-t-1}\\
>& (t+1)(k_1-t)(k_2-t+1)-\frac{(k_1-t)(k_2-t+1)}{2}-(k_1-t)(tk_2-t^2+1)-\frac{k_1-t}{2}\\
=&(k_1-t)(k_2+t-2)/2>0
\end{align*}}%
due to
{\small\begin{align*}
\frac{t(k_1-t)(k_1-t-1)}{n-t-1}\leq \frac{t(k_1-t)(k_1-t-1)}{(t+1)(k_2-t)(k_1-t+1)(k_2-t+1)}<\frac{k_1-t}{2}.
\end{align*}}%

Since $n\geq(t+1)(k_2-t)(k_1-t+1)(k_2-t+1)+t+1$ again, we have
{\small\begin{align*}
&\tilde{g}_6(k_1,k_2, n,t)-\tilde{g}_3(k_2,k_1, n,t)\\
>&\frac{n-t}{k_2-t}-\frac{(k_1-t)(k_2-t+1)}{2}-(k_2-t)(tk_1-t^2+1)-\frac{t(k_2-t)(k_2-t-1)}{n-t-1}\\
>&(t+1)(k_1-t+1)(k_2-t+1)-\frac{(k_1-t)(k_2-t)}{2}-\frac{k_1-t}{2}-(k_2-t)(tk_1-t^2+1)-1\\
>&(t+1)(k_1-t+1)(k_2-t)-\frac{(k_1-t)(k_2-t)}{2}-(k_2-t)(tk_1-t^2+1)\\
=&(k_2-t)(k_1+t)/2>0
\end{align*}}%
due to
{\small\begin{align*}
\frac{t(k_2-t)(k_2-t-1)}{n-t-1}\leq \frac{t(k_2-t)(k_2-t-1)}{(t+1)(k_2-t)(k_1-t+1)(k_2-t+1)}<1.
\end{align*}}%
Therefore, the required result follows.   $\qed$
\begin{lemma}\label{s-ineq-lem-8}
We have $g_1(k_1,k_2,n,t)>\max\{g_4(k_1,k_2,n,t),\ g_4(k_2,k_1,n,t)\}.$
\end{lemma}
\proof It suffices to prove that  $g_6(k_1,k_2,n,t)>\max\{g_4(k_1,k_2,n,t),\ g_4(k_2,k_1,n,t)\}$ by Lemma~\ref{s-ineq-lem-1}. $(k,\ell)=(k_1, k_2)$ or $(k_2,k_1)$. By (\ref{s-eq-4}), observe that
{\small\begin{align}\label{s-tg4}
\tilde{g}_4(k,\ell,n,t)=\frac{(k-t+1)^2(\ell-t-1)}{k-t}{t+2\choose 2}\left(1+\frac{1}{n-t-1}\right).
\end{align}}%

Since $k_1\geq k_2$, it is routine to check that $\frac{(k_1-t+1)(k_2-t-1)}{k_1-t}< k_2-t+1$, implying that
{\small\begin{align}\label{s-eq-g4-1}
\tilde{g}_4(k_1,k_2,n,t)<& (k_1-t+1)(k_2-t+1){t+2\choose 2}\left(1+\frac{1}{n-t-1}\right)\nonumber\\
<&(k_1-t+1)(k_2-t+1){t+2\choose 2}+1
\end{align}}%
from (\ref{s-tg4}) and $n\geq (t+1)^2(k_1-t+1)(k_2-t+1)+t+1$. By (\ref{s-tg6}) and (\ref{s-eq-g4-1}), observe that
{\small\begin{align*}
&\tilde{g}_6(k_1,k_2,n,t)-\tilde{g}_4(k_1,k_2,n,t)\\
>&n-t+\frac{n-t}{k_2-t}-\frac{(k_1-t)(k_2-t+1)}{2}-(k_1-t+1)(k_2-t+1){t+2\choose 2}-1\\
>&n-t-1-(k_1-t+1)(k_2-t+1)\cdot\frac{t^2+3t+3}{2}+\frac{n-t}{k_2-t}>0
\end{align*}}%
from $n\geq (t+1)^2(k_1-t+1)(k_2-t+1)+t+1$ and $(t+1)^2>(t^2+3t+3)/2$. Then $g_6(k_1,k_2,n,t)>g_4(k_1,k_2,n,t)$ holds.

Next, we show that $g_6(k_1,k_2,n,t)>g_4(k_2,k_1,n,t)$ holds. By (\ref{s-tg6}) and (\ref{s-tg4}), we have
{\small\begin{align*}
&\left(\tilde{g}_6(k_1,k_2,n,t)-\tilde{g}_4(k_2,k_1,n,t)\right)(k_2-t)(k_2-t+1)^{-1}\\
>&n-t-\frac{(k_1-t)(k_2-t)}{2}-(k_1-t-1)(k_2-t+1){t+2\choose 2}\left(1+\frac{1}{n-t-1}\right)\\
>&n-t-(k_1-t+1)(k_2-t+1)\left(\frac{t^2+3t+3}{2}+\frac{t^2+3t+2}{2(n-t-1)}\right)>0
\end{align*}}%
due to $n\geq (t+1)^2(k_1-t+1)(k_2-t+1)+t+1$ and $(t+1)^2>(t^2+3t+3)/2$. Hence, the required result follows. $\qed$
\begin{lemma}\label{s-ineq-lem-10}
We have $g_1(k_1,k_2,n,t)>g_5(k_1,k_2,n,t)$.
\end{lemma}
\proof From (\ref{s-eq-5}), observe that $\tilde{g}_5(k_1,k_2,n,t)=(t+1)^2(k_1-t+1)(k_2-t+1).$ Since $\max\{t+1,\ k_2-t\}\geq\frac{(k_2-t)(t+1)}{k_2-t+1}+\frac{k_2-t}{k_2-t+1},$ we have $n\geq \max\{t+1,\ k_2-t\}\cdot(t+1)(k_1-t+1)(k_2-t+1)+t+1\geq (t+1)^2(k_1-t+1)(k_2-t)+(t+1)(k_1-t+1)(k_2-t)+t+1$, implying that
{\small\begin{align*}
&(\tilde{g}_6(k_1,k_2,n,t)-\tilde{g}_5(k_1,k_2,n,t))(k_2-t)(k_2-t+1)^{-1}\\
>& n-t-\frac{(k_1-t)(k_2-t)}{2}-(t+1)^2(k_1-t+1)(k_2-t)>0
\end{align*}}%
due to (\ref{s-tg6}). Hence, the required result follows by Lemma~\ref{s-ineq-lem-1}. $\qed$

\section*{Acknowledgement}

M. Lu is supported by the National Natural Science Foundation of China (12171272), B. Lv is supported by National Natural Science Foundation of China (12071039, 12131011), K. Wang is supported by the National Key R\&D Program of China (No. 2020YFA0712900) and National Natural Science Foundation of China (12071039, 12131011).

%\end{CJK*}


\addcontentsline{toc}{chapter}{Bibliography}
\begin{thebibliography}{99}

\bibitem{a1}
N. Alon, L. Babai and H. Suzuki, Multilinear polynomials and Frankl-Ray-Chaudhuri-Wilson type intersection theorems, \emph{J. Combin. Theory Ser. A} 58 (1991) 165--180.
	
\bibitem{Ahlswede-Khachatrian-1996}
R. Ahlswede and L.H. Khachatrian, The complete nontrivial-intersection theorem for systems of finite sets, \emph{J. Combin. Theory Ser. A} 76 (1996) 121--138.

\bibitem{Ahlswede-Khachatrian-1997}
R. Ahlswede and L.H. Khachatrian, The complete intersection theorem for systems of finite sets, \emph{European J. Combin.} 18 (1997) 125--136.

\bibitem{Bey-2005} C. Bey, On cross-intersecting families of sets, \emph{Graphs Combin.} 21 (2005) 161--168.

\bibitem{Borg-2009} P. Borg, A short proof of a cross-intersection theorem of Hilton, \emph{Discrete Math.}
309 (2009), 4750--4753.

\bibitem{Borg-2014} P. Borg, The maximum sum and the maximum product of sizes of cross intersecting
families, \emph{European J. Combin.} 35 (2014), 117--130.

\bibitem{Borg-2016} P. Borg, The maximum product of weights of cross-intersecting families, \emph{J. London Math. Soc.} (2) 94 (2016) 993--1018.

\bibitem{Cao-set} M. Cao, B. Lv and K. Wang, The structure of large non-trivial t-intersecting families of finite sets. \emph{European J. Combin.} 97 (2021), 103373.

%\bibitem{Deza-Frankl-1983} M. Deza and P. Frankl, The Erd\H{o}s-Ko-Rado theorem--22 years later, \emph{SIAM J. Algebraic Discrete Methods} 4 (1983) 419--431.

\bibitem{Erdos-Ko-Rado-1961-313}
P. Erd\H{o}s,  C. Ko and R. Rado, Intersection theorems for systems of finite sets, \emph{Quart. J. Math. Oxf. Ser.} (2) 12(48) (1961) 313--320.

\bibitem{Frankl-1976} P. Frankl, On Sperner families satisfying an additional condition, \emph{J. Combin. Theory Ser. A} 20 (1976) 1--11.

%\bibitem{Frankl-1977} P. Frankl, Extremal set systems, PhD thesis, Hungarian Academy of Science, 1977 (in Hungarian).

\bibitem{Frankl-1978}
P. Frankl, The Erd\H{o}s-Ko-Rado theorem is true for $n = ckt$, in: Combinatorics, Vol. I, Proc. Fifth Hungarian Colloq., Keszthey, 1976, in: Colloq. Math. Soc. J\'{a}nos Bolyai, vol. 18, North-Holland, 1978, pp. 365--375.

\bibitem{Frankl-1978-1}
P. Frankl, On intersecting families of finite sets, \emph{J. Combin. Theory Ser. A} 24 (1978) 146--161.

%\bibitem{Frankl-1980}
%P. Frankl, On intersecting families of finite sets, \emph{Bull. Aust. Math. Soc.} 21 (1980) 363--372.

\bibitem{Frankl-1987}
P. Frankl, The shifting technique in extremal set theory, in: C. Whitehead (Ed.), Combinatorial Surveys, Cambridge Univ. Press, London, New York, 1987, pp. 81--110.

\bibitem{Frankl-Furedi-1986}
P. Frankl, Z. F\"{u}redi, Nontrivial intersecting families, \emph{J. Combin. Theory Ser. A} 41 (1986) 150--153.

\bibitem{Frankl--Furedi-1991}
P. Frankl and Z. F\"{u}redi, Beyond the Erd\H{o}s-Ko-Rado theorem, \emph{J. Combin. Theory Ser. A} 56 (1991) 182--194.

\bibitem{Frankl--Lee--Siggers--Tokushige-2014} P. Frankl, S. J. Lee, M. Siggers and N. Tokushige, An Erd\H{o}s-Ko-Rado theorem for cross $t$-intersecting families, \emph{J. Combin. Theory Ser. A} 128 (2014) 207--249.

%\bibitem{Frankl--Ota--Tokushige-1991}
%P. Frankl, K. Ota, N. Tokushige, Covers in uniform intersecting families and a counterexample to a conjecture of %lov\'{a}sz, \emph{J. Combin. Theory Ser. A} 74 (1996) 33--42.

\bibitem{Frankl--Tokushige-1992} P. Frankl, N. Tokushige, Some best possible inequalities concerning cross-intersecting families, \emph{J. Combin. Theory Ser. A} 61 (1992) 87--97.

\bibitem{Frankl-Tokushige-2002} P. Frankl, N. Tokushige, Weighted $3$-wise $2$-intersecting families, \emph{J. Combin. Theory Ser. A} 100 (2002) 94--115.

\bibitem{Frankl-Tokushige-2005} P. Frankl, N. Tokushige, Random walks and multiply intersecting families, \emph{J. Combin. Theory Ser. A} 109 (2005) 121--134.

\bibitem{Frankl--Tokushige-2011} P. Frankl, N. Tokushige, On $r$-cross intersecting families of sets, \emph{Combin. Probab. Comput.} 20 (2011) 749--752.

%\bibitem{Furedi-1988}
 %Z. F\"{u}redi, Matchings and covers in hypergraphs, \emph{Graphs Combin.} 4 (1988) 115--206.

%\bibitem{Furuya--Takatou-1988}
% M. Furuya, M. Takatou, Covers in $5$-uniform intersecting families with covering number three, \emph{Australas. J. %Combin.} 55 (2013) 249--262.

%\bibitem{GK}
%C. Godsil and M. Karen, Erd\H{o}s-Ko-Rado Theorems: Algebraic Approaches, Cambridge University Press, 2015.

\bibitem{Han-Kohayakawa}
J. Han and Y. Kohayakawa, The maximum size of a non-trivial intersecting uniform family that is not a subfamily of the Hilton-Milner family, \emph{Proc. Amer. Math. Soc.} 145(1) (2017) 73--87.

\bibitem{Hilton-1977} A.J.W. Hilton, An intersection theorem for a collection of families of subsets of a finite set, \emph{J. London
Math. Soc.} (2) 15 (1977) 369--376.

\bibitem{Hilton-Milner-1967}
A. Hilton and E. Milner, Some intersection theorems for systems of finite sets, \emph{Quart. J. Math. Oxford Ser.} (2) 18 (1967) 369--384.

%\bibitem{Katona-1964} G.O.H. Katona, Intersecting theorems for systems of finite sets, \emph{Acta Math. Acad. Sci. Hung.} 15 (1964) 329--337.

\bibitem{Kostochka-Mubayi}
A. Kostochka and D. Mubayi, The structure of large intersecting families, \emph{Proc. Amer. Math. Soc.} 145 (6) (2017) 2311-2321.

\bibitem{Liu-Zhang-Li-Zhang-2016} J. Liu, S. Zhang, S. Li, H. Zhang, Set systems with $k$-wise $L$-intersections and
codes with restricted Hamming distances, \emph{European J. Combin.} 58 (2016) 166--180.

\bibitem{Matsumoto-Tokushige-1989} M. Matsumoto and N. Tokushige, The exact bound in the Erd\H{o}s-Ko-Rado theorem for cross intersecting
families, \emph{J. Combin. Theory Ser. A} 52 (1989) 90--97.

\bibitem{Moura-1999}  L. Moura, Maximal $s$-wise $t$-intersecting families of sets: kernels, generating sets, and enumeration. \emph{J. Combin. Theory Ser. A} 87 (1999) 52--73.

\bibitem{O-V-2021} J. O'Neill and J. Verstra\"{e}te, Non-trivial $d$-wise intersecting family, \emph{J. Combin. Theory Ser. A} 178 (2021) 105369.

\bibitem{Pyber-1986} L. Pyber, A new generalization of the Erd\H{o}s-Ko-Rado theorem, \emph{J. Combin. Theory Ser. A} 43 (1986)
85--90.

\bibitem{Tokushige-2006} N. Tokushige, Extending the Erd\H{o}s-Ko-Rado theorem, \emph{J. Combin. Des.} 14 (2006) 52--55.

\bibitem{Tokushige-2007} N. Tokushige, The maximum size of $3$-wise $t$-intersecting families, \emph{European J. Combin.} 28 (2007) 152--166.

\bibitem{Tokushige-2007-2}  N. Tokushige, EKR type inequalities for $4$-wise intersecting families, \emph{J. Combin. Theory Ser. A} 114 (2007) 575--596.

\bibitem{Tokushige-2010} N. Tokushige, On cross $t$-intersecting families of sets, \emph{J. Combin. Theory Ser. A} 117 (2010) 1167--1177.

\bibitem{Tokushige-2013} N. Tokushige, The eigenvalue method for cross $t$-intersecting families, \emph{J. Algebr. Comb.} 38 (2013)
653--662.

\bibitem{Wang-Zhang-2011} J. Wang and H. Zhang, Cross-intersecting families and primitivity of symmetric
systems, \emph{J. Combin. Theory Ser. A} 118 (2011), 455--462.


\bibitem{Wilson-1984}
R.M. Wilson, The exact bound in the Erd\H{o}s-Ko-Rado theorem, \emph{Combinatorica} 4 (1984) 247--257.
\end{thebibliography}
\end{document}